\documentclass [12pt]{amsart}

\usepackage{amssymb,amsxtra,amsfonts}
\usepackage{epsfig}             %
\usepackage{graphics}

\usepackage[colorlinks]{hyperref}

\usepackage{verbatim}  
\usepackage{bbding}   

\newenvironment{ppb}[1]
{\ \!\!\!\!\!\!\!\!\!\!\!\!\!\!\!\!\!\!\!\!\!\!\!\!\!\!\!\!\!\!\!\!\!\!\!\!\!\!\!\! {\bf PPB------------------------------------------------------------------------------------------------PPB}\newline \tiny {#1}
\  \newline\normalsize\phantom{f}\!\!\!\!\!\!\!\!\!\!\!\!\!\!\!\!\!\!\!\!\!\!\!\!\!\!\!\!\!\!\!\!\!\!\!\!\!\!\!\! {\bf PPB------------------------------------------------------------------------------------------------PPB}\newline}{}

\def\reE@DeclareMathSymbol#1#2#3#4{%
    \let#1=\undefined
    \DeclareMathSymbol{#1}{#2}{#3}{#4}}
\DeclareSymbolFont{symbolsC}{U}{txsyc}{m}{n}
\SetSymbolFont{symbolsC}{bold}{U}{txsyc}{bx}{n}
\DeclareFontSubstitution{U}{txsyc}{m}{n}
\reE@DeclareMathSymbol{\strictiff}{\mathrel}{symbolsC}{76}

\newcommand\beq{\begin{equation}}
\newcommand\eeq{\end{equation}}
\newcommand\bal{\begin{align*}}
\newcommand\eal{\end{align*}}   
\newcommand\bmx{\left(\begin{matrix}}
\newcommand\emx{\end{matrix}\right)}
\newcommand\bsmx{\left(\begin{smallmatrix}}
\newcommand\esmx{\end{smallmatrix}\right)}

\newcommand{\spq}{/\!\!/}

\providecommand{\spqa}[1]{\underset{#1}{/\!\!/}}

\newcommand{\st}{\ \bigl\vert\ }

\providecommand{\Rad}{\text{\rm Rad}}

\def\part#1{\frac{\partial\phantom{q}}{\partial#1}}

\newcommand {\flb}{\lbrack\!\lbrack}
\newcommand {\frb}{\rbrack\!\rbrack}
\newcommand {\flp}{(\!(}
\newcommand {\frp}{)\!)}

\newcommand{\glue}[1]{\underset{#1}{\strictiff}}
\newcommand{\fus}{\circledast}
\newcommand{\fusion}[1]{\underset{#1}{\circledast}}

\newcommand{\MB}{\mathcal{M}_{\text{\rm B}}}

\newcommand{\Lie}{{\mathop{\rm Lie}}}

\DeclareMathOperator{\ISto}{{\IS}to} 

\newcommand{\pap}[2]{{\ _{#1}\cA_{#2}}}

\newcommand{\papk}[3]{  \ _{#1}^{\phantom{#3}}\cA_{#2}^{#3}     }
\newcommand{\gah}{\pap{G}{H}}  
\newcommand{\gat}{\pap{G}{T}}  
\newcommand{\gaho}{\papk{G}{H}{1}}   
\newcommand{\gato}{\papk{G}{T}{1}}   
\newcommand{\gaht}{\papk{G}{H}{2}} 

\newcommand{\gatt}{\papk{G}{T}{2}} 

\newcommand{\gatr}{\papk{G}{T}{r}} 
\newcommand{\gahr}{\papk{G}{H}{r}}

\DeclareMathOperator{\pr}{pr}

\newcommand{\Prod}{\prod}

\DeclareMathOperator{\Hom}{Hom}         
\DeclareMathOperator{\Aut}{\mathop{\rm Aut}}

\newcommand{\GL}{{\mathop{\rm GL}}}

\newcommand{\U}{{\rm {U}}}	
\newcommand{\SO}{{\mathop{\rm SO}}}

\newcommand{\irr}{{\rm irr}}

\DeclareMathOperator{\Rep}{\rm Rep}
\DeclareMathOperator{\End}{End}

\newcommand{\diagonals}{{\mathop{\rm diagonals}}}

\newcommand{\reg}{{\mathop{\footnotesize\rm reg}}}

\newcommand{\hk}{{hyperk\"ahler }}   

\newcommand{\ba}{{\bf a}}

\newcommand{\bH}{{\bf H}}

\newcommand{\bQ}{{\bf Q}}
\newcommand{\bs}{{\bf S}}
\newcommand{\bS}{{\bf S}}

\newcommand{\IA}{\mathbb{A}}
\newcommand{\IB}{\mathbb{B}}
\newcommand{\IC}{\mathbb{C}}
\newcommand{\ID}{\mathbb{D}}

\newcommand{\IG}{\mathbb{G}}
\newcommand{\IH}{\mathbb{H}}

\newcommand{\IP}{\mathbb{P}}

\newcommand{\IS}{\mathbb{S}}

\newcommand{\IZ}{\mathbb{Z}}

\newcommand{\cA}{\mathcal{A}}
\newcommand{\cB}{\mathcal{B}}
\newcommand{\cC}{\mathcal{C}}

\newcommand{\cK}{\mathcal{K}}

\newcommand{\cM}{\mathcal{M}}

\newcommand{\cO}{\mathcal{O}}

\newcommand{\R}{\mathcal{R}}
\newcommand{\cR}{\mathcal{R}}

\newcommand{\gM}{       \mathfrak{M}     }

\newcommand{\g}{       \mathfrak{g}     }

\newcommand{\lt}{\mathfrak{t}}

\newcommand{\gl}{       \mathfrak{gl}     } 

\renewcommand{\k}{\mathfrak{k}}

\newcommand{\wt}{\widetilde}

\newcommand{\wh}{\widehat}

\newcommand{\al}{\alpha}

\newcommand{\de}{\delta}
\newcommand{\De}{\Delta}

\newcommand{\Ga}{\Gamma}

\newcommand{\la}{\lambda}
\newcommand{\La}{\Lambda}

\newcommand{\Si}{\Sigma}

\renewcommand{\bar}{\overline}

\makeatletter
 \newlength{\typesize}
\makeatother

\newlength{\vvoff}
\newlength{\hhoff}

\def\mapright#1{\smash{
        \mathop{\longrightarrow}\limits^{#1}}}

\def\underset#1#2{\ \smash{\mathop{ #2 }\limits_{#1}}\ }

\newcommand{\pf}{\begin{bpf}}

\newcommand{\pfms}{\begin{bpfms}}
\newcommand{\epf}{\end{bpf}\hfill$\square$\\}           
\newcommand{\epfms}{\end{bpfms}\hfill$\square$\\}       

\newcommand{\idea}{\begin{bidea}}

\newcommand{\eidea}{\end{bidea}\hfill$\square$\\}           

\newcommand{\sk}{\begin{bsk}}    

\newcommand{\esk}{\end{bsk}\hfill$\square$\\}           
\newcommand{\sketch}{\begin{bsketch}}

\newcommand{\esketch}{\end{bsketch}\hfill$\square$\\}

\newtheorem {hypo}{\bf\hspace{-\parindent}Hypothesis}
\newtheorem {thm}[hypo]{Theorem}   
\newtheorem {prop}[hypo]{Proposition}

\newtheorem {cor}[hypo]{Corollary}
\newtheorem {lem}[hypo]{Lemma}

\newtheorem {defn}[hypo]{Definition}

\theoremstyle{remark}\newtheorem{rmk}[hypo]{Remark}

\renewenvironment{ppb}[1]{}{}

\begin{document}

\title{Poisson varieties from Riemann surfaces}
\author{Philip Boalch}%

\begin{abstract}
Short survey based on talk  at the Poisson 2012 conference.
The main aim is to describe and give some examples of wild character varieties 
(naturally generalising the character varieties of Riemann surfaces by allowing
more complicated behaviour at the boundary), 
their Poisson/symplectic structures 
(generalising both the Atiyah--Bott approach and the quasi-Hamiltonian approach), and the wild mapping class groups.
\ppb{
Some links to other parts of mathematics will also be sketched, including:
1) The fact that standard quantum group quantizes a simple wild character variety, 
2) hyperk\"ahler geometry of wild character varieties and the corresponding extension of nonabelian Hodge theory,
3) the new (Poisson) discrete group actions that arise upon enriching the notion of ``smooth complex algebraic curve'' curve into an ``irregular curve'', and the resulting generalisation of the mapping class group.}
\end{abstract}

\maketitle

\section{Introduction}

To fix ideas recall that there is already a much-studied class of Poisson varieties attached to Riemann surfaces, as follows.

Let $G$ be a connected complex reductive group, such as $G=\GL_n(\IC)$,  with a  nondegenerate invariant symmetric bilinear form on its Lie 
algebra.
Then there is a Poisson variety, the $G$-character variety, 
attached to any Riemann surface $\Si$:
\beq\label{eq: curves to spaces}
\Si \qquad \mapsto \qquad \Hom(\pi_1(\Si),G)/G,
\eeq
taking the space of conjugacy classes of representations in $G$ of the fundamental group.
In brief these spaces are complex symplectic if $\Si$ is compact (without boundary) and in general 
they are Poisson with symplectic leaves obtained by fixing the conjugacy class of the monodromy around each boundary component.

The first symplectic approach to such spaces of fundamental group representations was analytic
\cite{AB83} and subsequently there were many alternative, more algebraic, approaches
to these character varieties, such as  
\cite{Gol84, Kar92,FR,  huebschmann-duke95, 
Jeffrey-95, Weinstein-floermem, AMR, GHJW, AMM}.
Several applications are discussed in the surveys \cite{Ati90, Aud95long}.

One of the interesting features is that these character varieties are finite dimensional Poisson manifolds which, unlike the vast majority of  the Poisson manifolds appearing in classical mechanics, cannot be constructed out of finite dimensional cotangent bundles or coadjoint orbits.
Rather, they are finite dimensional spaces  naturally obtained from infinite dimensions, 
and this is the approach taken by Atiyah--Bott \cite{AB83} 
(and this fact provided motivation for the above listed authors to 
seek a purely finite-dimensional approach).

However if one reads a little further in the Poisson literature it soon becomes apparent that there are other important Poisson manifolds which do not fit into this framework 
(of either being constructed out of finite dimensional cotangent bundles/coadjoint orbits or from moduli spaces of flat connections).

For example the Drinfeld--Jimbo quantum group $U_q(\g)$ 
has an integral structure in which we can put $q=1$ 
and obtain the following Poisson manifold: 
\beq\label{eq: plgp}
G^* := \{ (b_+,b_-)\in B_+\times B_- \st \de(b_+)\de(b_-)=1 \}
\eeq
where $B_\pm\subset G$ are a pair of opposite Borel subgroups and $\de:B_\pm\to T$ is the natural  projection onto the maximal torus $T=B_+\cap B_-$.
(Thus if $G=\GL_n(\IC)$ we could take $T$ to be the diagonal subgroup and $B_\pm$ to be the upper/lower triangular subgroups and then $\de$ is the map taking the diagonal part.)
This statement is proved in \cite{DKP},\cite{DP1565}  Theorem p.86 \S12.1, and 
an earlier version of this result at the level of formal groups is due to Drinfeld, cf. \cite{Drin86} \S3.
These Poisson manifolds are also discussed in \cite{STS83, FRT88, AlekMalk94, YKS97}.

Two facts suggest that there might be some link between $G^*$ and spaces of connections:

1)  $G^*$ is constructed in terms of the algebraic groups $T,B_\pm$ (rather than say their Lie algebras) and algebraic groups were first considered, autonomously from Lie theory, as ``Galois groups'' for differential equations (see e.g. Chevalley's MR reviews of \cite{kolchin46a,kolchin46b,kolchin48}), so naively one might guess that any natural space  involving algebraic groups will have a `differential equations' counterpart\footnote{Another instance of this philosophy is the fact that the Grothendieck simultaneous resolution is a moduli space of (``logahoric'') connections \cite{logahoric, hdr}.},

2) The symplectic leaves of $G^*$ 
are obtained by fixing the conjugacy class of the product $b_-^{-1}b_+\in G$ (which is reminiscent of fixing the conjugacy class of monodromy around the boundary above).

Such considerations lead to the following:

\noindent{\bf Question:} Is there a class of Poisson varieties associated to connections on Riemann surfaces, generalising the above spaces of fundamental group representations, and which includes the spaces $G^*$?

The main aim of this note is to describe the answer to this question 
(from \cite{thesis, smid, smapg, bafi,saqh02}), and some more recent developments \cite{fission, gbs}.
In brief the answer is to consider more general algebraic connections than those parameterised by representations of the fundamental group---this algebraic viewpoint means that such moduli spaces are often overlooked in differential geometry, even though they parameterise objects well-known 
in algebraic geometry and the theory of differential equations, 
as will be explained.
It turns out that spaces of such (more general) connections also have natural Poisson/symplectic structures.
The classical ideas of Riemann and others on Fuchsian differential equations and the behaviour of their solutions (monodromy) lead to the character varieties, whereas 
more recent ideas, understanding the behaviour of solutions of non-Fuchsian equations, lead to the wild character varieties we will describe here.

The original motivation (\cite{thesis, smid}) for this work was to better understand the theory of {\em isomonodromic deformations}, a class of nonlinear differential equations.
Our secondary aim in this note is to explain how this comes about, i.e. how certain nonlinear differential equations may be understood geometrically from this set-up, of associating Poisson moduli spaces to Riemann surfaces.

For example, in the original context of the character varieties 
sketched above, suppose 
we vary the surface $\Si$ in a smooth family over a base $\IB$.
Then, via the map \eqref{eq: curves to spaces}, we get a family of moduli spaces parameterised by $\IB$, and it turns out they fit together to form the fibres of a 
fibre bundle over  $\IB$, which has a natural flat connection on it, i.e. a flat Ehresmann/nonlinear connection 
(the fibres, the character varieties, are not linear spaces).
Integrating this flat connection yields a Poisson action of the fundamental group of the base on the fibres, i.e. a homomorphism
$$\pi_1(\IB,b)\to \Aut_{\text{Poisson}}\bigl(\Hom(\pi_1(\Si_b),G)/G\bigr).$$
Typically $\pi_1(\IB,b)$ will be a braid or mapping class group 
(the fundamental group of a moduli space of pointed Riemann surfaces).
Our second aim is to explain how to generalise this story: 
the notion of surface with marked points should (and will) be generalised and this leads
 to new spaces of deformation parameters whose fundamental groups 
act on the generalised character varieties by Poisson automorphisms. 
The simplest example involves the $G$-braid group, and gives a geometric derivation of the (so-called) quantum Weyl group action on $G^*$.

\subsection{From compact groups to \hk manifolds}

In fact most of the references \cite{AB83, Gol84, Kar92,FR,  huebschmann-duke95, Jeffrey-95, Weinstein-floermem, AMR, GHJW, AMM}   deal with character varieties for compact groups, i.e. 
 spaces of the form  $\Hom(\pi_1(\Si),K)/K$
for compact Lie groups $K$ such as $K=U_n$, rather than the complex reductive groups $G=K_\IC$.
Much of the motivation for the case of compact groups comes from the theorem of Narasimhan--Seshadri \cite{NarSes}
relating such spaces to
moduli spaces of stable vector bundles on  $\Si$:

\begin{thm}
Suppose $\Si$ is compact and  $K=U_n$. Then the space 
$$\Hom^{\irr}(\pi_1(\Si),K)/K$$
of irreducible representations of the fundamental group of $\Si$ is a K\"ahler manifold 
and is isomorphic to the moduli space of stable degree zero rank $n$ vector bundles on $\Si$. 
\end{thm}
In particular $\Hom^{\irr}(\pi_1(\Si),K)/K$ has an underlying (real) symplectic structure and the results cited above gave various different constructions of this symplectic form.
Whereas the K\"ahler metric depends on the complex structure of the surface $\Si$ the symplectic structure is topological, and depends only on the underlying real two-manifold.

Here we are more interested in the case of complex reductive groups (i.e. groups which are obtained by complexifying compact groups). 
In this case, by complexifying the above constructions, one obtains a complex/holomorphic  symplectic form on $\Hom^{\irr}(\pi_1(\Si),G)/G$, which is still topological.
The corresponding extension of the theorem of Narasimhan--Seshadri 
was proved by 
is due to Hitchin, Donaldson, Corlette and  Simpson 
(\cite{Hit-sde, Don87, Cor88, Sim-hbls}):

\begin{thm}\label{thm: hdcs}
Suppose $\Si$ is compact and  $G=\GL_n(\IC)$. 
Then the space $$\Hom^{\irr}(\pi_1(\Si),G)/G$$
of irreducible representations of the fundamental group of $\Si$ is a \hk manifold 
and is isomorphic to the moduli space of stable degree zero rank $n$ Higgs bundles on $\Si$. 
\end{thm}

In brief, the complexified Atiyah--Bott construction realises 
$\Hom^{\irr}(\pi_1(\Si),G)/G$ as a complex symplectic quotient of an infinite dimensional vector space,
and thereby explains why it is symplectic.
Hitchin's approach to Theorem \ref{thm: hdcs} is a strengthening of this symplectic quotient construction, into a \hk quotient.
Since our focus is the Poisson/symplectic geometry, this will not be discussed further here, 
although it should be noted that one of the 
original motivations in \cite{smid}
for generalising the Atiyah--Bott symplectic structure was to obtain new \hk manifolds, 
by also extending Theorem \ref{thm: hdcs}:
this was done in \cite{wnabh}---i.e. the generalised/irregular Atiyah--Bott symplectic quotient (of \cite{thesis,smid}) was strengthened there into a \hk quotient---this \hk story is surveyed in \cite{ihptalk}.

\section{Riemann--Hilbert}

Suppose $G=\GL_n(\IC)$. To understand how to generalise 
$\Hom(\pi_1(\Si),G)/G$
it is worth carefully considering the following:

\noindent{\bf Question:} What objects are parameterised by $\Hom(\pi_1(\Si),G)/G$?

First, a differential geometer may well quote:

\begin{thm}
Consider the set of pairs $(V,\nabla)$ where $V\to \Si$ is a $C^\infty$ rank $n$ complex vector bundle 
and $\nabla$ is a flat connection on $V$. Then 
there is a bijection
$$\{ (V,\nabla) \}/\text{isomorphism} \qquad \cong \qquad  \Hom(\pi_1(\Si),G)/G$$
obtained by taking a connection to its monodromy  representation.
\end{thm}

Secondly,  a complex geometer may well note that the $(0,1)$ part of such a flat connection determines a holomorphic structure on the bundle $V$, and in that holomorphic structure the connection is holomorphic, and thus derive the stronger result:

\begin{thm}
Consider the set of pairs $(V,\nabla)$ where $V\to \Si$ is a rank $n$ holomorphic vector bundle 
and $\nabla$ is a holomorphic connection on $V$. Then 
there is a bijection
$$\{ (V,\nabla) \}/\text{isomorphism} \qquad \cong \qquad  \Hom(\pi_1(\Si),G)/G$$
obtained by taking a connection to its monodromy  representation.
\end{thm}

Thirdly an expert on differential equations may well recall that the notion of monodromy was introduced by Riemann in his work on the hypergeometric differential equation, and one would like to use such monodromy data to classify {\em algebraic} differential equations (or connections).
For example Hilbert's 21st problem is often interpreted in the following way:

Fix  distinct complex numbers $a_1,\ldots,a_m\in \IC$ 
and consider the space of connections
\beq\label{eq: lifted log conns}
\wt\cM^* \ :=\  \Bigl\{ d-\sum_1^m \frac{A_i}{z-a_i} dz \st A_i\in \End(\IC^n) \Bigr\}
\ \cong\  \End(\IC^n)^{m}
\eeq
on the trivial rank $n$ bundle on the $m$-punctured plane $\Si=\IC\setminus\{a_i\}$.
Choosing a basepoint $x\in \Si$ and taking the monodromy of the connections around loops based at $x$ 
then yields a holomorphic map
\beq\label{eq: RH}
\End(\IC^n)^{m}\cong \wt \cM^*\quad  \mapright{\text{RH}} \quad \Hom(\pi_1(\Si,x),G) \cong \GL(\IC^n)^{m}
\eeq
between two spaces of the same dimension.
A key point is that by considering such a space of algebraic connections we can talk about the algebraic structure of such moduli spaces (i.e. the algebraic structure on the space of coefficients $A_i$)
and thus see that the Riemann--Hilbert map is a natural holomorphic map between two different algebraic varieties 
(generalising the matrix exponential when $m=1$).
Since the spaces on both sides are of the same dimension, and the map is bijective on an open neighbourhood of the trivial connection, one would like to know:

Question: Is this Riemann--Hilbert map surjective?

\noindent
or even:

Question: is there a way to make this into a precise, bijective correspondence?

This first question was answered in the negative (for $n,m$ sufficiently large) 
by Bolibruch (see \cite{Bol94, beauville-RH}).
For the second question, the expert on differential equations, or an algebraic geometer, may well now recall Deligne's Riemann--Hilbert correspondence:

\begin{thm}[\cite{Del70}]\label{thm: dRH}
Suppose $\Si$ is a smooth complex algebraic curve (not necessarily compact).
Consider the set of pairs $(V,\nabla)$ where $V\to \Si$ is a rank $n$ complex algebraic vector bundle 
and $\nabla$ is an algebraic connection on $V$ which has regular singularities at infinity. 
Then there is a bijection
$$\{ (V,\nabla) \}/\text{isomorphism} \qquad \cong \qquad  \Hom(\pi_1(\Si),G)/G$$
obtained by taking a connection to its monodromy  representation.
\end{thm}

Thus Deligne tells us exactly what are the ``algebraic differential equations'' (and the precise equivalence relation amongst them), that are 
parameterised  by the fundamental group representations. (In fact he establishes a more general result in arbitrary dimensions.)
Notice that, unlike in the holomorphic context, not all algebraic connections are considered: the condition of regular singularities at infinity means that the horizontal sections grow at most polynomially at each puncture (or equivalently that the bundle has an extension across any puncture in which the connection has at most a first order pole).
In fact in the genus zero case Plemelj had proved earlier that all representations arise as monodromy of a regular singular connection (cf. \cite{beauville-RH}).

To clarify how vast is the class of connections which are missing in this correspondence, note the following:

\begin{prop}[cf. \cite{BV83} p.28]
Suppose $\nabla$ is an algebraic connection on an algebraic vector bundle $V\to \Si$
which locally at a puncture takes the form
$$d-\left(\frac{A_k}{z^k} + \frac{A_{k-1}}{z^{k-1}}+\cdots\right)dz$$
where $k\ge 2$ and $A_k$ is not nilpotent.
Then $\nabla$ does not have regular singularities.
\end{prop}

This suggests, therefore, where to look to find a generalisation of the character varieties: we should look at moduli spaces of connections which do not have regular singularities, i.e. irregular connections, and describe such moduli spaces in terms of monodromy-type data.

Fortunately a Riemann--Hilbert correspondence (on smooth algebraic curves for $G=\GL_n(\IC)$) for such connections has been worked out  (starting with the preliminary work of Birkhoff \cite{birkhoff-1913}, and leading to e.g. \cite{sibuya77, jurkat78, BJL79, malg-book, DMR-ci}), although it is more complicated than
the regular singular case (one should bear in mind that much of the theory of Riemann surfaces that we now take for granted, such as the uniformization theorem, was originally worked out in the context of understanding multivalued  solutions of linear differential equations).

Before discussing spaces of generalised monodromy data 
in \S\ref{sn: wcv}, we will first review the tame case in the next section.
Section \ref{sn: sympl strs} will then describe 
the purely algebraic approach to the symplectic/Poisson structures on the 
wild character varieties and \S\ref{sn: ff} will discuss some of the ideas in the proof.
Many examples are then described in
\S\ref{sn: egs}.

\begin{rmk}
The above theorems imply that any irregular algebraic connection is {\em holomorphically} isomorphic to a regular singular algebraic connection.
But this kills all the ``new'' moduli that occur in the irregular case: 
rather the aim is to construct interesting new moduli spaces 
generalising the character varieties, show they have all the key properties of the usual ``tame'' case 
and derive new features as well 
(such as the irregular braid group actions generalising the usual mapping class group).
In particular the holomorphic maps \eqref{eq: RH} admit suitable irregular generalisations.
The examples below show that the ``new'' moduli behave just like the
familiar $\pi_1$ representations, and the simplest wild character varieties  are actually {\em isomorphic} to tame cases (but there are many genuinely new moduli spaces, 
and the familiar spaces of fundamental group representations are just the tip of the iceberg). 
\end{rmk}

\section{Quasi-Hamiltonian approach to tame character varieties}

Now we will recall the quasi-Hamiltonian approach to tame character varieties, which gives a clean algebraic approach to the  Atiyah--Bott symplectic structures. 

Suppose $\Si$ is a compact Riemann surface of genus $g$ with $m\ge 1$ boundary components $\partial_1,\ldots,\partial_m$, and let $G$ be 
a complex reductive group as above (such as $G=\GL_n(\IC)$).
Let $\MB(\Si):=\Hom(\pi_1(\Si),G)/G$ denote the corresponding character variety.
Choose conjugacy classes $\cC_1,\ldots,\cC_m\subset G$ and let 
$$\MB(\Si,\cC):=\Hom_\cC(\pi_1(\Si),G)/G \ \subset\  
\MB(\Si)$$
denote the corresponding symplectic leaf of the character variety (parameterising connections with monodromy around $\partial_i$ in $\cC_i$). The notation $\MB$ alludes to the fact that these are the ``Betti moduli spaces'' in the terminology of nonabelian cohomology.

Upon choosing suitable loops generating the fundamental group, 
the character variety $\Hom_\cC(\pi_1(\Si),G)/G$ clearly ``looks like'' a multiplicative symplectic quotient:

$$\Hom_\cC(\pi_1(\Si),G)/G \cong \qquad\qquad\qquad\qquad\qquad\qquad\qquad\qquad\qquad\qquad\qquad\qquad$$
\nopagebreak
$$\qquad\qquad\{ (A_i,B_i,M_j)\in G^{2g}\times \Prod\cC_j\st 
[A_1,B_1]\cdots[A_g,B_g] M_1\cdots M_m = 1\}/G$$

\noindent
where $G$ acts by diagonal conjugation, and $[a,b]=aba^{-1}b^{-1}$.

The quasi-Hamiltonian theory
 of Alekseev--Malkin--Meinrenken \cite{AMM} sets up a formalism of ``Lie group valued moment maps'' such that this actually 
{\em is} a multiplicative symplectic quotient\footnote{
In fact \cite{AMM} only considered compact groups and some of their proofs do not work directly in the noncompact case---the case of complex reductive groups was first considered 
in \cite{saqh02}, 
and a new approach to quasi-Hamiltonian geometry via Dirac structures, which works equally well in the 
noncompact case, was 
worked out later (see  \cite{ABM-purespinors} and references therein).}---and the symplectic structure on 
the character variety may be obtained algebraically in this way, from the structure of ``quasi-Hamiltonian $G$-space''
(upstairs) on $G^{2g}\times \Prod\cC_j$.
As mentioned in \cite{AMM} this approach is a reinterpretation of an earlier construction of \cite{Jef94, huebschmann-duke95, GHJW} (and the natural 
quasi-Hamiltonian two-form on conjugacy classes already played 
an important role in \cite{GHJW}).

Without going too much into the details (see \cite{AMM}) they showed that 
the space
$$\ID(G) = G\times G$$
(the internally fused double)
is a  ``quasi-Hamiltonian $G$-space'' with moment map
$\mu(A,B) = [A,B]\in G$,
as is any conjugacy class $\cC\subset G$, with moment map given by the inclusion 
(a multiplicative analogue of a coadjoint orbit $\cO\subset \g^*\cong \g$).
They also defined a ``fusion'' operation
which implies the product
$$G^{2g}\times \Prod\cC_j = \ID(G)^g\times \cC_1\times\cdots\times\cC_m$$
inherits the structure of quasi-Hamiltonian $G$-space with moment map given by the product of the moment maps on each component, i.e.
$$\mu(A_i,B_i,M_j)=[A_1,B_1]\cdots[A_g,B_g] M_1\cdots M_m\in G$$
precisely as required.
Finally, they defined a reduction operation, a multiplicative analogue of
the symplectic quotient, which implies the character variety,
now identified with $\mu^{-1}(1)/G$,  inherits a genuine symplectic structure.
In symbols:
\beq\label{eq: qhisom}
\Hom_\cC(\pi_1(\Si),G)/G\ \cong\  
\ID^{\fus{}g}\fus\cC_1\fus\cdots\fus\cC_m\spq G
\eeq
where $\fus$ is the fusion product, and $\spq$ denotes the 
multiplicative symplectic quotient.

Thus by developing the theory of group valued moment maps it is possible to show the
character varieties are finite dimensional {multiplicative} symplectic quotients\footnote{unfortunately to obtain the full \hk structure, so far, it seems one has to start in infinite dimensions.}.

The isomorphism \eqref{eq: qhisom}
depends on a choice of generators of the fundamental group, so does not 
reveal the intrinsic nature of the quasi-Hamiltonian set-up.
A more intrinsic statement is as follows.

Choose a basepoint $b_i$ on the $i$th boundary component $\partial_i$ 
for each $i=1,\ldots, m$.
Let $$\Pi=\Pi_1(\Si,\{b_1,\ldots,b_m\})$$
be the fundamental groupoid of $\Si$ with basepoints $\{b_i\}$, i.e. the groupoid of homotopy classes of paths in $\Si$ whose endpoints are in the set of chosen basepoints.

\begin{thm}\label{thm: amm gpoid spaces}
The space $\Hom(\Pi,G)$ of homomorphisms from the groupoid $\Pi$ to the group $G$ is a smooth affine variety and is a quasi-Hamiltonian $G^m$-space, with moment map
$$\mu: \Hom(\Pi,G) \to G^m$$
having $i$th component the monodromy around $\partial_i$ based at $b_i$.
\end{thm}

This statement is the starting point of  \cite{gbs} and  is essentially just a rephrasing of (the complexification of) 
the main result of \cite{AMM}.
It immediately implies that the quotient 
$$\Hom(\Pi,G)/G^m,$$
which is isomorphic to 
$\Hom(\pi_1(\Si),G)/G$, 
is a Poisson variety\footnote{Similarly 
Theorem \ref{thm: amm gpoid spaces} also immediately implies that
$\Hom(\Pi,G)/G^{m-1} \cong \Hom(\pi_1(\Si,b_1),G)$ is a quasi-Poisson $G$-manifold 
(cf. \cite{AKM}).
This quasi-Poisson structure is mapping class group invariant, since the quasi-Hamiltonian structure on  
$\Hom(\Pi,G)$ is mapping class group invariant (either by direct computation as in \cite{AMM}, or since it comes from the intrinsic Atiyah--Bott approach). 
A new approach to this invariance in the $\GL_n$ case 
appeared recently  in \cite{MassTur12}.}.
And moreover it follows directly that 
the symplectic leaves are the spaces $\mu^{-1}(\cC)/G^m$ for any conjugacy class
$\cC=(\cC_1\ldots,\cC_m)\subset G^m$, i.e. they are the spaces 
$\MB(\Si,\cC)=\Hom_\cC(\pi_1(\Si),G)/G$ with fixed local monodromy conjugacy classes.
It is this statement that we will generalise.

\section{Wild character varieties}\label{sn: wcv}

In this section we will define the wild character varieties, which generalise the character varieties 
$\MB(\Si,\cC)\subset \MB(\Si)$ considered above (which will henceforth be referred to as the ``tame character varieties'').
In the case $G=\GL_n(\IC)$, much of this section is essentially a recollection of known results on the irregular Riemann--Hilbert correspondence. 
We have rephrased them in a convenient way, extended them to general complex reductive groups and shown how to obtain algebraic varieties (parameterising generalised monodromy data). 

\ppb{
When \cite{AMM} first appeared on the arxiv in 1997, the author had recently extended the Atiyah--Bott symplectic structure to give a similar analytic description 
of many moduli spaces of irregular connections, as an infinite dimensional symplectic quotient. 
These moduli spaces also have a finite-dimensional description, again in a form that looks just like a multiplicative symplectic quotient. 
This raised the question to extend the quasi-Hamiltonian approach so that these more general spaces {\em are} multiplicative symplectic quotients.

In this section we will review the finite-dimensional description of these spaces:}

In brief one enriches the fundamental group representation by adding in the 
Stokes data. There are various ways to present this extra data and we will describe an approach involving local systems on an associated real surface (see the references in \cite{smid, bafi, gbs} 
such as \cite{BJL79, malg-book, MR91, L-R94} for more background).

The reader might skip forward %
 to see some simple examples of the new quasi-Hamiltonian spaces that appear, before mastering the details of the general definitions here.

\subsection{Irregular curves.}

Fix a connected complex reductive group $G$, such as $G=\GL_n(\IC)$, and a maximal torus $T\subset G$.
A wild character variety $\MB(\Si)$ 
is then associated to an object called an ``irregular curve'' (or ``wild Riemann surface'') $\Si$. 
This generalises the notion of curve with marked points, in order to encode some of the boundary conditions for irregular connections.
An irregular curve $\Si$ (for fixed $G$) is a triple $(\Si,\ba,\bQ)$ consisting of a smooth compact complex algebraic curve $\Si$, plus some marked points 
$\ba=(a_1,\ldots,a_m)\subset \Si$, for some $m\ge 1$, and an ``irregular type'' $Q_i$ at each marked point $a_i$.
The extra data   is the irregular type: an irregular type at a point 
$a\in \Si$ is an arbitrary element 
$Q\in \lt(\wh \cK_a)/\lt(\wh \cO_a)$ where $\lt=\Lie(T)$.
If we choose a local coordinate $z$ vanishing at $a$ then
$\wh \cK_a=\IC\flp z \frp, \wh \cO_a=\IC\flb z \frb$ so that 
$Q\in \lt\flp z \frp/\lt\flb z \frb$
and we may write
$$Q = \frac{A_r}{z^r}+\cdots + \frac{A_1}{z}$$ 
for  some elements $A_i\in \lt$, for some integer $r\ge 1$.
Note that we place no restriction on the coefficients $A_i$.
(The tame/regular singular case is the case when each 
$Q_i=0$, and in this case the wild character variety coincides with the usual space of fundamental group representations.)
There are several variations of this definition 
(bare irregular curves, twisted irregular curves) discussed in \cite{gbs} but they will not be needed here.
Although it may seem very strange at first, this ``irregular curve'' viewpoint is very convenient, since it turns out that the extra moduli in the irregular types behaves just like the moduli of the curve,
and so it makes sense to include it with the choice of the curve with marked points right from the start. 

\subsection{Connections on irregular curves.}
Given an irregular curve $\Si=(\Si,\ba,\bQ)$ we may consider algebraic connections on algebraic $G$-bundles on $\Si^\circ = \Si\setminus\{\ba\}$
such that near each puncture $a_i$ there is a local trivialisation such that the connection takes the form 
\beq\label{eq: fnf}
dQ_i + \La(z)dz/z,\eeq
for some $\g$ valued map $\La(z)$ (holomorphic across the puncture).
In other words we consider connections whose irregular part is determined up to isomorphism by $dQ_i$.
From now on, given an irregular curve $\Si$, the words ``a connection on $\Si^\circ$'' will always mean such an algebraic connection on a principal $G$-bundle on $\Si^\circ$. 
Note that local solutions to such connections near $a_i$ involve the essentially singular term $\exp(Q_i)$.

\subsection{Irregular Riemann--Hilbert and Stokes local systems}
One way to state the irregular Riemann--Hilbert correspondence in this context is that: 
Given an irregular curve $\Si$, the category (groupoid) of connections on  $\Si^\circ$ is equivalent to the category 
of {\em Stokes $G$-local systems} determined by $\Si$ (cf. \cite{gbs} Appendix A).
This statement generalises %
Deligne's equivalence (Theorem \ref{thm: dRH})  between regular singular connections on algebraic vector bundles on 
$\Si^\circ$, and local systems on $\Si^\circ$.
Thus Stokes local systems are topological objects that classify such algebraic connections.

A Stokes $G$-local system is defined as follows. Given the irregular 
curve $\Si$ we may define two real surfaces with boundary 
$$\wt\Si\hookrightarrow \wh\Sigma \twoheadrightarrow  \Si$$ as follows:
$\wh\Si$ is the real oriented blow-up of $\Si$ at the points of $\ba$, so that each point $a_i$ is replaced by a circle $\partial_i$ parameterising the real oriented direction emanating from $a_i$, and the boundary of  $\wh \Si$ is $\partial\wh \Si=\partial_1\sqcup\cdots \sqcup \partial_m$.
In turn each irregular type $Q_i$ canonically determines three pieces of data 
(see Appendix \ref{apx: sd and sg} below or \cite{gbs} for full details):

1) a connected reductive group $H_i\subset G$, the centraliser of $Q_i$ in $G$,

2) a finite set $\IA_i\subset \partial_i$ of singular/anti-Stokes directions for all $i=1,\ldots,m$, and

3) a unipotent group $\ISto_d\subset G$ normalised by $H_i$, the Stokes group of $d$, for each $d\in \IA_i$ for all $i=1,\ldots,m$.

The surface $\wt \Si\subset \wh \Si$ is defined by puncturing $\wh \Si$ once at a point $e(d)$ in its interior near each singular direction $d\in \IA_i\subset \partial_i$, for all $i=1,\ldots,m$.
For example we could fix a small tubular neighbourhood of $\partial_i$ (a ``halo'')  $\IH_i\subset \wh \Si$  (so $\IH_i$ an annulus), and 
choose the extra puncture $e(d)\in\wh \Si$ to be on the {\em interior} boundary of $\IH_i$ near the singular direction $d$ (as pictured in Figure \ref{fig: halo}).

\begin{figure}[ht]
	\centering
	\input{halo.on.curve2.pstex_t}
	\caption{The surface $\wt \Si$ with halo drawn}\label{fig: halo}
\end{figure}

A Stokes $G$-local system for the irregular curve $\Si$ consists of 
 a $G$-local system on $\wt \Si$, with a flat reduction to $H_i$ in $\IH_i$ for each $i=1,\ldots,m$, such that the local monodromy around $e(d)$ is in 
$\ISto_d$ for any basepoint in $\IH_i$, for all $d\in \IA_i$ for all $i$.

Thus, via the irregular Riemann--Hilbert correspondence, the classification of connections with given irregular types is reduced to the classification of Stokes $G$-local systems.
Thus it remains to describe the set of isomorphism classes of Stokes local systems.
This goes as follows.
Choose a basepoint $b_i\in \IH_i$ for each $i$ and let 
$\Pi$ be the fundamental groupoid of $\wt\Si$ with basepoints 
$\{b_1,\ldots,b_m\}$.
Then we may consider the space $\Hom(\Pi,G)$ of $G$-representations of $\Pi$ 
and the subspace $\Hom_\IS(\Pi,G)$ of ``Stokes representations'' which are the representations such that the local monodromy around $\partial_i$ (based at $b_i$) is in $H_i$ and the local monodromy around $e(d)$ (based at $b_i$) is in $\ISto_d$,  for all $d\in\IA_i$, for all $i$.
The group $\bH:=H_1\times\cdots\times H_m$ acts naturally on $\Hom_\IS(\Pi,G)$.
The basic classification statement is then:

\begin{thm}
The set of isomorphism classes of Stokes $G$-local systems for the irregular curve $\Si$ (and thus the set of isomorphism classes of connections on $\Si^\circ$) is naturally in bijection with the set of $\bH$ orbits in $\Hom_\IS(\Pi,G)$.
\end{thm}

A slightly different approach to defining $\wt \Si$ will be described in Appendix \ref{apx: morecan}.

\subsection{Wild character varieties.}
By choosing paths generating $\Pi$ it is easy to see that the space $\Hom_\IS(\Pi,G)$ of Stokes representations  is naturally a smooth affine variety.
The  wild character variety $\MB(\Si)$ is defined to be the affine quotient
$$ \MB(\Si) := \Hom_\IS(\Pi,G)/\bH$$
taking the variety associated to the ring of $\bH$-invariant functions.
This means that the points of $\MB(\Si)$ correspond bijectively to the {\em closed} $\bH$ orbits in  $\Hom_\IS(\Pi,G)$.
Beware that in this notation $\Si$ now denotes an {\em irregular} curve.

Note that if $\Si$ is tame, i.e. if all the irregular types $Q_i$ are zero, then $\bH=G^m$ and the corresponding character variety
$$ \MB(\Si) = \Hom_\IS(\Pi,G)/ G^m \cong \Hom(\pi_1(\Si^\circ),G)/G$$
coincides with the usual space of conjugacy classes of  representations of the fundamental group of the punctured curve.

In turn suppose we choose a conjugacy class $\cC\subset \bH$ in the group $\bH$, so that $\cC = (\cC_1,\ldots,\cC_m)$ where $\cC_i\subset H_i$ is a conjugacy class in $H_i$.
Then we may consider the (locally closed) subvariety 
$$\MB(\Si,\cC) \subset \MB(\Si)$$
defined by restricting to Stokes representations whose monodromy around $\partial_i$ (based at $b_i$) is in the conjugacy class $\cC_i\subset H_i$, for $i=1,\ldots,m$.
Again, in the tame case this specialises to the usual character variety with fixed conjugacy classes of local monodromy (since then $\cC\subset \bH=G^m$).

\begin{rmk}
Beware that in the wild case in general the Stokes representations in $\MB(\Si,\cC)$ do {\em not} correspond to algebraic connections with fixed local monodromy conjugacy classes around the punctures. Rather, it is the conjugacy class of the so-called ``formal monodromy'' (in $H_i$) that we are fixing, rather than the actual local monodromy (in $G$) of the connections around the punctures.
These two notions coincide in the tame case.
\end{rmk}

\section{Symplectic and Poisson structures}\label{sn: sympl strs}

Thus in summary the previous section defined the wild character varieties: 
 an algebraic variety $\MB(\Si)$ attached to an irregular curve $\Si$, and for any choice of conjugacy class
$\cC\subset \bH$, a subvariety $\MB(\Si,\cC)\subset\MB(\Si)$.
If we choose suitable paths generating the fundamental groupoid $\Pi$ considered above, then again these spaces ``look just like''  a multiplicative symplectic quotient
(cf. \eqref{eq: big fusion prod redn} below).
The author noticed this after reading about the
specific explicit examples considered 
by Jimbo--Miwa--Ueno in \cite{JMU81} 
(see especially their equation (2.46)),
although any mention of symplectic/Poisson geometry 
was absent there.

Thus, after \cite{AMM} appeared, it was natural to ask if one could extend the theory of Lie group valued moment maps to incorporate such more general spaces (and thus give a finite dimensional/algebraic approach to the irregular Atiyah--Bott symplectic structures of \cite{thesis, smid, wnabh}).
This is indeed the case:

\begin{thm}(\cite{saqh02,saqh,fission,gbs})\label{thm: qh stokes reps}
Given an irregular curve $\Si$, the space $\Hom_\IS(\Pi,G)$
of Stokes representations is a smooth affine variety, which has the structure of quasi-Hamiltonian $\bH$-space with moment map 
$$\mu:\Hom_\IS(\Pi,G) \to \bH;\qquad 
\rho\mapsto \{\rho(\partial_i) \st i=1,\ldots, m\}$$
taking the formal monodromies.
\end{thm}

The general quasi-Hamiltonian/quasi-Poisson yoga then immediately implies that the ring of $\bH$-invariant functions on 
$\Hom_\IS(\Pi,G)$ is a Poisson algebra (i.e. that the affine quotient $\MB(\Si)$ is a Poisson variety), with symplectic leaves given by the subvarieties $\MB(\Si,\cC)$ with fixed formal monodromy conjugacy classes.
In other words the spaces $\MB(\Si,\cC)$
are the multiplicative symplectic quotients of $\Hom_\IS(\Pi,G)$
by $\bH$ (at values of the moment map determined by $\cC$).
This statement is a natural generalisation of Theorem 
\ref{thm: amm gpoid spaces} to the irregular case: 
if all the irregular types are zero, 
Theorem \ref{thm: qh stokes reps} specialises to Theorem \ref{thm: amm gpoid spaces}.

This result was obtained in 2002
\cite{saqh02, saqh} in the case where the leading term $A_r\in\lt$ of each irregular type was {\em regular semisimple}, i.e. off of all of the root hyperplanes---this condition implies (but is not equivalent to) the centralizer groups $H_i$ each being a maximal torus.
Note that in such cases (whenever $\bH$ is abelian) 
Theorem \ref{thm: qh stokes reps} implies that $\Hom_\IS(\Pi,G)$ is itself a complex symplectic manifold (the quasi-Hamiltonian two-form is then symplectic).
These symplectic structures (on $\Hom_\IS(\Pi,G)$ and $\MB(\Si,\cC)$) 
were constructed analytically earlier in \cite{thesis, smid}, and this analytic construction extends to the general case 
(cf. \cite{wnabh} for the \hk strengthening of it, to the general case, for general linear groups).
As we will see below, many interesting examples (such as $G^*$) appear already in this special case (with $A_r$ regular).

\section{Fusion and Fission}\label{sn: ff}

Here we will mention some of the ideas in the proof of Theorem \ref{thm: qh stokes reps}.
One of the key ideas of \cite{AMM} (closely related to the Hamiltonian loop group viewpoint 
\cite{Don-bvym, MW-jdg})
is to construct moduli spaces of flat connections on an arbitrary surface by fusion, i.e. 
gluing together pairs of pants, 
thereby enabling an induction with respect to the genus and the number of poles.
In the irregular case, the local picture at a pole is still too complicated to construct 
moduli spaces in one go: a new operation appears, enabling a further induction with respect to the order of the pole (or more precisely with respect to the number of ``levels''---cf.  Appendix \ref{apx: sd and sg}).

\subsection{Fusion}

\begin{figure}[h]
\begin{center}
\includegraphics[angle=0, width=11cm, trim=0 0 0 0, clip]{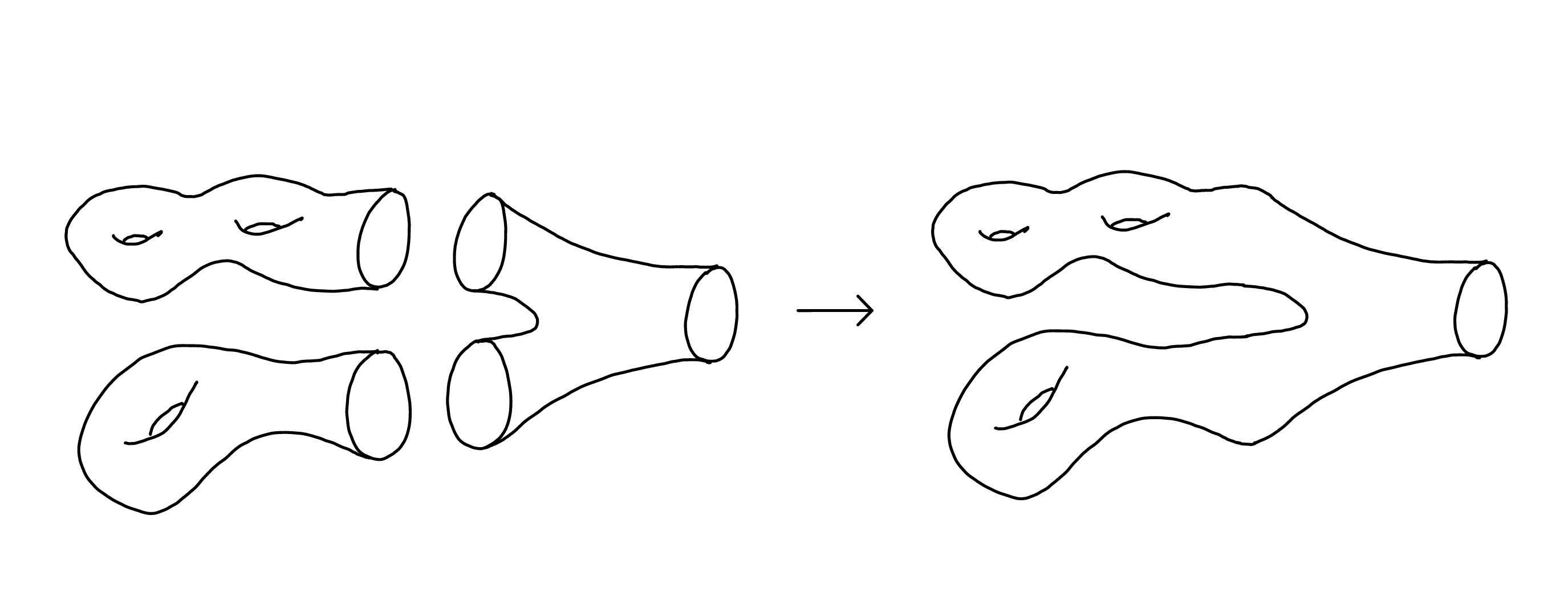}
\caption{Fusion}\label{fig: fusion}
\end{center}
\end{figure}

The fusion idea (see Figure \ref{fig: fusion}) 
gives a way to chop the Riemann surface up in to three-holed spheres (pairs of pants), and the quasi-Hamiltonian incarnation of this idea enables the full 
symplectic moduli space to be constructed by gluing together simple pieces:

\begin{thm}(\cite{AMM})
All of the tame character varieties may be constructed 
inductively
by gluing together conjugacy classes and pairs of pants.
\end{thm}

\subsection{Fission}

Using the fusion idea in the presence of irregular singular connections still enables an induction with respect to the genus and the number of poles, allowing a reduction to the case of just one pole (on a disc).
However, this case is still too complicated to handle directly and another idea is needed: A new operation ``fission'' enables an induction with respect to the {order} of the poles, as mentioned above. 

Choose an integer $r\ge 1$ and an element $A\in \lt$ in the Cartan subalgebra.
Let $\Delta$ be the unit disc in the complex plane with coordinate $z$ and let $Q_0=A/z^r$ be an irregular type at $z=0$ with just one term. 
Consider the (local) irregular curve 
$\Si= (\De,0,Q_0)$.

Let $\gahr=\Hom_\IS(\Pi(\Si),G)$ 
be the space of Stokes representations attached to $\Si$.
As a special case of Theorem \ref{thm: qh stokes reps}
$\gahr$ is a quasi-Hamiltonian $G\times H$ space, 
where $H=C_G(A)$ is the centralizer of $A$.
The spaces $\gahr$  are called ``fission spaces''.

\begin{thm}(\cite{gbs})\label{thm: fission}
All of the wild character varieties may be constructed 
inductively
by gluing together conjugacy classes, pairs of pants and fission spaces.
\end{thm}

If we have an arbitrary irregular type 
$Q=\sum A_i/z^i$ (still on a disk), the key point is 
to prove that
the corresponding space
\beq\label{eq: one pole}
\cA(Q)\ :=\ \Hom_\IS(\Pi(\De,0,Q),G) \ \cong\  
G\times H\times\!\! \Prod_{d\in \IA(Q)} \ISto_d(Q)
\eeq
of Stokes representations
is a quasi-Hamiltonian $G\times H$-space, where 
$H= C_G(Q)$).
This follows  since $\cA(Q)$
may be obtained by gluing a sequence of fission spaces end to end, for a nested sequence of centralizer groups $H$ 
(see \eqref{eq: nesting} below).

Thus Theorem \ref{thm: fission} enables us to reduce the general case of Theorem \ref{thm: qh stokes reps} to the special case of the fission spaces $\gahr$.
In turn this is much easier since the spaces $\gahr$ are quite simple:

\begin{lem}\label{lem: easier}(see \cite{gbs}, proof of Lemma 7.11, and \cite{bafi})
There are a pair of opposite parabolic subgroups 
$P_\pm\subset G$ with $H=P_+\cap P_-$ and an isomorphism
$$\gahr \quad \cong \quad G\times (U_+\times U_-)^r\times H$$
where $U_\pm = \Rad_u(P_\pm)$ are the unipotent radicals of the parabolics.
\end{lem}

Thus the main remaining step is to establish:

\begin{thm}
For any $r\ge 1$ and pair of opposite parabolics 
$P_\pm=H\cdot U_\pm\subset G$, the space
$$G\times (U_+\times U_-)^r\times H$$
is a quasi-Hamiltonian $G\times H$ space.
\end{thm}

This was established in \cite{saqh02, saqh} 
in the case when $H=T$ was a maximal torus,
in \cite{fission} in the 
case $r=1$ (for any $H$), and in \cite{gbs} in general.

In these coordinates the moment map $\mu=(\mu_G,\mu_H):\gahr\to G\times H$
is %
$$\mu_G(C,\bs,h) = C^{-1} h S_{2r}\cdots S_2 S_1 C$$
$$\mu_H(C,\bs,h) = h^{-1}$$
where $C\in G, h\in H$ and 
$\bS\in (U_+\times U_-)^{r}$ with $\bS = (S_1,\ldots,S_{2r})$ and
$S_{\text{even}}\in U_-$ and $S_{\text{odd}}\in U_+.$
The  action of $G\times H$ on $\gahr$ given by 
$(g,k)(C,\bS,h)  = (kCg^{-1}, k\bS k^{-1}, khk^{-1})$
where $(g,k)\in G\times H$ and $k\bS k^{-1} = (kS_1k^{-1},\ldots,kS_{2r}k^{-1})$.

Further \cite{saqh02} explains how these quasi-Hamiltonian structures may be {\em derived} 
from the irregular Atiyah--Bott approach of \cite{thesis, smid}, and this explains why the resulting quasi-Hamiltonian structures have lots of nice properties 
(which may be checked algebraically, as in \cite{gbs}).
Note that after \cite{thesis, smid} appeared Woodhouse \cite{Woodh-sat} and Krichever
\cite{Krich-imds}
wrote down some explicit formulae for these ``irregular Atiyah--Bott'' symplectic structures
on Stokes data (in the case where the leading term $A_r$ at each pole was regular semisimple, for general linear groups).

One may see that the fission spaces yield new operations on the category of quasi-Hamiltonian spaces, 
as  follows.
First, the usual fusion picture may be phrased as follows. 
Let $S$ denote the three-holed sphere. 
This yields a quasi-Hamiltonian $G^3$-space 
$M(S) = \Hom(\Pi_1(S),G)$, where $\Pi_1(S)$ denotes the fundamental groupoid of $S$ with one basepoint on each boundary component. %
Similarly if $\Si_1,\Si_2$ are two surfaces with just one boundary component then the spaces 
$M(\Si_1), M(\Si_2)$
are quasi-Hamiltonian $G$-spaces.
Then fusion amounts to the following {\em gluing}:
$$M(\Si_1)\fus M(\Si_2) = M(\Si_1)\glue{G} M(S) \glue{G} M(\Si_2)$$

\noindent
where the symbol $\glue{}$ denotes the gluing---which on the level of surfaces just glues the corresponding boundary circles together (cf. \cite{fission} \S5).  
Since $M(S)$ is a quasi-Hamiltonian $G^3$-space, and each gluing absorbs a factor of $G$, the result is a quasi-Hamiltonian $G$-space, as expected.

Now, for the fission spaces, typically $H$ will factor as a product of groups (e.g. if $G=\GL_n(\IC)$ then $H$ is a ``block diagonal'' subgroup). 
Suppose $H=H_1\times H_2$ for definiteness (the generalisation to arbitrarily many factors is immediate).
Thus $\gahr$ is a quasi-Hamiltonian $G\times H_1\times H_2$ space.
For example this enables us to construct a quasi-Hamiltonian $G$-space 
$$M_1\glue{H_1} \gahr \glue{H_2} M_2$$

\noindent
for any integer $r$, out of quasi-Hamiltonian $H_i$-spaces $M_i$ ($i=1,2$), e.g. we could take  $M_i = \Hom(\pi_1(\Si_i),H_i)$.
Thus the fission spaces yield many new operations on the category of quasi-Hamiltonian spaces (without fixing the group $G$ beforehand).
One may picture these operations as indicated in Figure \ref{fig: Y}.

\begin{figure}[ht]
	\centering
	\input{annulus.and.Y.pstex_t}
	\caption{Fission%
}\label{fig: Y}
\end{figure}

First (as above)
the fission spaces may be pictured in terms of surfaces as on the left in 
Figure \ref{fig: Y}.
Alternatively, perhaps more accurately, if the centraliser group $H= H_1\times H_2$ is a product, then, once we cross the boundary of the halo and break the group from $G$ to $H$, the two factors of $H$ are completely independent and the inner annulus may break into two independent sheets, one for each factor, as indicated on the right of 
 Figure \ref{fig: Y} (cf. \cite{fission, gbs}).
This $Y$-shaped piece defines a new operation enabling one to glue various surfaces (with $H_1,H_2$ and $G$ connections on them, respectively).

In particular the spaces $\cA(Q)$ attached to 
an arbitrary irregular type 
$Q=\sum_1^r A_i/z^i$ (on a disk, as in \eqref{eq: one pole}), may be obtained by gluing fission spaces end to end, as follows:
Consider the nested  chain of connected complex reductive groups 
$$H=H_1\subset H_2\subset \cdots H_r \subset G$$
where $H_i\subset G$ is the centraliser of $\{A_r,\ldots,A_i\}\subset \g$.
Then \cite{gbs} Prop 7.12 shows there is an explicit ``nesting'' isomorphism of quasi-Hamiltonian $G\times H$-spaces
\beq\label{eq: nesting}
\cA(Q)\quad\cong \quad \cA(r)\glue{H_r}\cdots \glue{H_3}\cA(2)\glue{H_2} \cA(1)
\eeq
where $\cA(i) =\!\! \papk{H_{i+1}}{H_{i}}{i}$ is a fission space
for the groups $H_i\subset H_{i+1}$.
The dimension of $\cA(Q)$ may be computed as follows:
$$\dim(\cA(Q)) = \dim(G)+\dim(T) + 
\#\{\al\in \cR\st \deg (\al\circ Q)=0\} + \sum_{\al\in \cR}\deg(\al\circ Q)$$
where $\cR\subset \lt^*$ is the set of roots of $G$ and the ``degree''
of $\al\circ Q(z) \in \IC[z^{-1}]$ is its pole order, i.e. the degree of the polynomial $\al\circ Q(z^{-1})$. 
In the tame  case ($Q=0$) the space $\cA(Q)$ 
is just the double $D(G)$ of \cite{AMM}, and if $Q$ has regular semisimple leading term $A_r$ then 
$\cA(Q)$ is the space $\wt \cC/L$ of 
\cite{saqh02, saqh}, and is isomorphic to the fission space $\!\!\gatr$ and thus, via Lemma 
\ref{lem: easier}, to 
$G\times T\times (U_+\times U_-)^r$ 
with $U_\pm$ the unipotent radicals of a pair of opposite Borel subgroups.

Putting all this together,
given an arbitrary irregular curve $\Si$, 
 upon choosing suitable generators of the corresponding groupoid $\Pi=\Pi_1(\wt \Si)$
the resulting explicit description of the space  
$\Hom_\IS(\Pi,G)$ of Stokes representations is as follows:

\begin{align}\label{eq: big fusion prod}
\Hom_\IS(\Pi,G) \ & \cong\ \left(
\ID^{\fus g}\fusion{G}\cA(Q_1)
\fusion{G}\cdots\fusion{G} \cA(Q_m)\right)\spq{} G \\ & \nonumber \\
\ &\cong \ \left(\ID^{\fus g}\fusion{G}\cA(Q_1)
\fusion{G}\cdots\fusion{G} \cA(Q_{m-1})\right)\glue{G} \cA(Q_m)\nonumber
\end{align}

\noindent
where $g$ is the genus of $\Si$ (cf. \cite{gbs} (36)).
In particular (ignoring degenerate cases):
$$\dim\Hom_\IS(\Pi,G) = (2g-2)\dim(G) +\sum_1^m \dim\cA(Q_i).$$

Further, given a choice of conjugacy classes $\breve\cC_i\subset H_i$, i.e. a conjugacy class 
$\breve\cC=(\breve\cC_1,\ldots,\breve\cC_m)\subset \bH$ then
we can define a quasi-Hamiltonian $G$-space
\beq\label{eq: generalised cc}
\cC_i \ :=\  
\cA(Q_i)\glue{H_i}  \breve\cC_i  \ = \ 
\left(\cA(Q_i)\fusion{H_i}  \breve\cC_i\right)\spq H_i
\eeq
for each $i=1,\ldots,m$.
Thus if the singularity at $a_i$ is tame ($Q_i=0$) 
then $\cC_i= \breve\cC_i$ is a conjugacy 
class in $G=H_i$, but in general it is different, and can be highly 
non-trivial even if $\breve\cC_i$ is a point.
Then the symplectic wild character variety $\MB(\Si,\breve\cC)$ arises as the multiplicative symplectic reduction:
\beq\label{eq: big fusion prod redn}
\MB(\Si,\breve\cC)\ \cong\ \left(
\ID^{\fus g}\fusion{G}\cC_1
\fusion{G}\cdots\fusion{G} \cC_m\right)\spq{} G
\ \cong\ 
\Hom_\IS(\Pi,G)\spqa{\breve \cC'} \bH
 \eeq

\noindent
where $\breve \cC'\subset \bH$ is the conjugacy class inverse to $\breve\cC$.
These are the symplectic leaves of the Poisson wild character variety 
$\MB(\Si)=\Hom_\IS(\Pi,G)/ \bH$.

\ppb{
As above,
the fission spaces may be pictured in terms of surfaces as on the left in 
Figure \ref{fig: Y}.
Alternatively, perhaps more accurately, the centraliser group $H$ will typically be a product of groups, for example $H= H_1\times H_2$, and then, once we cross the boundary of the halo and break the group from $G$ to $H$, the two factors of $H$ are completely independent and the inner annulus may break into two independent sheets, one for each factor, as indicated on the right of 
 Figure \ref{fig: Y} (cf. \cite{fission, gbs}).
This $Y$-shaped piece defines a new operation enabling one to glue various surfaces (with $H_1,H_2$ and $G$ connections on them, respectively).

\begin{figure}[ht]
	\centering
	\input{annulus.and.Y.pstex_t}
	\caption{Fission%
}\label{fig: Y}
\end{figure}
}

\subsection{Fission varieties}

Note that the converse of Theorem \ref{thm: fission}
is not true: if we consider the class of symplectic/Poisson or quasi-Hamiltonian varieties which arise 
by gluing together conjugacy classes, pairs of pants and fission spaces, then one obtains a larger class of spaces than the 
wild character varieties (an example is described in \cite{fission}). 
We call this class of varieties the ``fission varieties'' (it is natural to include also the 
tame fission spaces of \cite{logahoric}).
This is in contrast to the tame case: by gluing together conjugacy classes and pairs of pants, only tame character varieties will be obtained.

Thus, in achieving our goal to construct the wild character varieties, a new building block was constructed, which leads to many new symplectic varieties, beyond those we were trying to construct.
It is not clear if all the symplectic manifolds in this larger class enjoy all of the good properties that the wild character varieties have (such as being \hk\!\!---cf. the conjecture in \cite{fission}).

\subsection{Weighted conjugacy classes}
The algebraic connections considered here are in fact much closer to considering meromorphic connections on $G$-bundles on the compact curve, i.e. allowing connections with poles 
(than to considering holomorphic connections on the open curve).
An example of this extended viewpoint is the bijective 
irregular Riemann--Hilbert statement of \cite{smid} Cor. 4.9.
Moreover in the full \hk story one should consider such extensions, as well as the 
more general story with compatible parabolic structures at the poles 
(cf. \cite{wnabh, ihptalk}), and in turn parahoric structures in 
the case of general reductive groups \cite{logahoric}.
We have skipped this here to simplify the presentation, but in terms of the Betti moduli spaces 
it amounts simply to replacing each of the conjugacy classes $\breve \cC_i\subset H_i$ above
by a ``weighted conjugacy class'' $\wh \cC_i$ for $H_i$, as in \cite{logahoric} \S4.
The key point is that $\wh \cC_i$ is still a quasi-Hamiltonian 
$H_i$-space, and so can be used to replace $\breve \cC_i$ in 
\eqref{eq: generalised cc}.
(Such spaces appear in the Brieskorn--Grothendieck--Springer simultaneous resolution, and such ideas lead to the fact that the resolution map is a quasi-Poisson moment map \cite{logahoric, hdr}.)
In the tame case for $\GL_n$ this amounts to adding a filtration preserved by the local monodromy, as was considered by Levelt \cite{levelt61} (2.2) and Simpson \cite{Sim-hboncc} for 
logarithmic connections on (parabolic) vector bundles,
enriching the notion of ``local system'' into 
the notion of ``filtered local system'' 
(in general it amounts to considering ``filtered Stokes $G$-local systems'' as in \cite{gbs}).

\section{Examples}\label{sn: egs}

First we will describe some examples where the centraliser group $H$ is a maximal torus (\cite{smid, saqh02}).
Suppose $B_\pm\subset G$ are a pair of opposite Borel subgroups, with unipotent radicals $U_\pm\subset G$.
Thus if $G=\GL_n(\IC)$ 
we could take $B_+$ to be the upper triangular subgroup of $G$
so that $U_+\subset B_+$ is the upper triangular subgroup
with $1$'s on the diagonal and $U_-\subset B_-$ are the corresponding lower triangular subgroups, so that $T=B_+\cap B_-$ is the diagonal torus.

Note that the product map
$U_-\times T \times U_+ \to G$ is an isomorphism onto its image, and this image 
$G^\circ=U_- T U_+ \subset G$ is a dense open subset of the group $G$.

\begin{cor}
The space $G^\circ$ is a Poisson manifold.
\end{cor}
\pf
It is a manifold since it is an open subset of $G$,
and it inherits a Poisson structure 
since it is isomorphic to the quotient 
$\gaho/G$ (with $r=1$ and $H=T$).
Indeed, in general $\gahr/G$ is a quasi-Poisson $H$-space, which means it is Poisson in the present case when $H$ is abelian. 
\epf

In other words, combining the above statements: the open subset 
$G^\circ$ inherits a natural Poisson structure from the irregular Atiyah--Bott symplectic structure.

It also follows that the symplectic leaves are obtained by fixing the conjugacy class of the monodromy around the outer boundary of the disc, i.e. the symplectic leaves are the 
intersections  $\cC\cap G^\circ$, for conjugacy classes $\cC\subset G$.

Now recall (from \eqref{eq: plgp}) 
the Poisson manifold $G^*\subset B_+\times B_-$ underlying the
 Drinfeld--Jimbo quantum group, and let $\wt G^*$ denote its universal cover.
Clearly the space $G^*$ is isomorphic to $U_+\times T\times U_-$
and so $\wt G ^* \cong U_+\times \lt \times U_-$.
Thus there are covering maps 
$$\wt{G}^*\ \to\  G^*\  \to\  G^\circ$$
taking $(u_+,\La,u_-)\in \wt{G}^*$
to $(su_+,s^{-1}u_-)\in G^*$  and to 
$u_-^{-1}tu_+\in G^\circ$
where $s=\exp(\pi i \La), t=s^2$.
Thus, since the coverings are isomorphisms at the level of tangent spaces, 
the Poisson structure on $G^o$ yields a Poisson structure on both 
$\wt G^*$ and on $G^*$.
(In fact this universal covering, including $\La$ and not just $t=\exp(2\pi i \La)$, was 
considered moduli theoretically from the start in \cite{smid}---this discrete choice of logarithm 
of the formal monodromy $t$
 corresponds to choosing an extension of the bundle across the puncture at $z=0$.)

This gauge theoretic Poisson structure on $G^*$ (coming from the fact it is a wild character variety) coincides  
with the Poisson structure on $G^*$ coming from the fact it is the quasi-classical limit of a quantum group:

\begin{thm}(\cite{smapg, bafi})
The natural Poisson structure on  $G^*$ arises from the 
geometry of irregular connections on a disk having poles of order two, and so the Drinfeld--Jimbo quantum group is the quantisation of a space of irregular connections.
\end{thm}

Similarly there are analogous spaces with more unipotent factors: 
The manifold $$ T\times(U_+\times U_-)^r$$
has a natural Poisson structure for any $r\ge 1$, since it is isomorphic to the quotient $\gatr/G$. 
The symplectic leaves are obtained by fixing the $G$-conjugacy class of the product $tS_{2r}\cdots S_2S_1\in G$, where $S_{even}\in U_-$, $S_{odd}\in U_+$, and $t\in T$.

The above examples arise from considering connections on a disc with just one marked point (and regular semisimple leading term of the irregular type).
If we consider two such discs and glue them along their boundary to obtain a Riemann sphere with two marked points, then one obtains a quasi-Hamiltonian $T\times T$ space, which is in particular a holomorphic symplectic manifold
(since $T\times T$ is abelian). 
In the case $r=1$ this yields (\cite{saqh02, saqh}) 
the Lu--Weinstein ``double symplectic groupoid'' \cite{LuW-sdgpoid},
and so gives a moduli theoretic realisation of a class of symplectic manifolds first constructed abstractly 
for different reasons. 
In this context the irregular Riemann--Hilbert map yields a transcendental map from the total space of the cotangent bundle of $G$, to the double symplectic groupoid (see \cite{saqh02, saqh}).

Alternatively we could glue on a disc with no marked points, thus yielding a quasi-Hamiltonian $T$-space 
(and so in particular a symplectic manifold).
This corresponds to considering the symplectic leaf in the Poisson manifold $\gatr/G$
determined by the identity conjugacy class in $G$. 
For example if we take $r=2$ 
then one obtains the following symplectic manifolds.

\begin{cor} (\cite{thesis, smid})
The intersection 
\beq\label{eq: oppintn}
  (U_+ U_-)\  \cap \ (U_- T U_+) \ \subset \ G
  \eeq
is a symplectic manifold, of dimension equal to the number of roots of $G$.
\end{cor}
\pf
It is a manifold since it is an open subset of $U_+ U_-\cong U_+\times U_-$.
Now consider the reduction $\cB := \gatt\spq{}G$ of the corresponding fission space by the action of $G$, at the value $1$ of the $G$ moment map.
This reduction is a quasi-Hamiltonian $T$-space and so symplectic, and easily 
identified with \eqref{eq: oppintn}.
\epf

The more recent developments \cite{fission, gbs} involve the fission spaces $\!\gahr$ for arbitrary subgroups $H\subset G$ that arise as Levi factors of parabolics $P\subset G$: in particular $H$ is again always a connected complex reductive group.
This, with $H$ nonabelian, allows us to glue several fission spaces end to end 
for sequences of reductive groups 
$H_1\subset H_2\subset G$ (and varying parameters $r$).

For example suppose we restrict to general linear groups (so in turn each Levi factor $H$ is a product of general linear groups).
Then the resulting fission spaces can be used \cite{cmqv} to set up
a precise multiplicative version of the theory of Nakajima quiver varieties.
This gives a new way to think about a large class of the wild character varieties as ``multiplicative versions'' of some quite well-known symplectic varieties.

In brief, given a graph $\Ga$ with nodes $I$ and an $I$-graded vector space $V=\bigoplus_I V_i$
define the vector space $\Rep(\Ga,V)$ of representations of the graph $\Ga$ on $V$ to be the space
of linear maps 
$$\Rep(\Ga,V) = \bigoplus _{e\in \bar{\Ga}}\Hom(V_{t(e)},V_{h(e)})$$
along each edge in each direction.
Here $\bar{\Ga}$ is the double of $\Ga$, i.e. the set of all possible oriented edges; given an oriented edge $e$, then $h(e),t(e)\in I$ are the nodes at its head and tail resp. 
Given a choice of an orientation of $\Ga$ then $\Rep(\Ga,V)$ inherits a natural symplectic structure and it has a Hamiltonian action of the group $\IG=\Prod_I \GL(V_i)$.
Given a choice of a parameter $\la\in \IC^I$, i.e. a scalar $\la_i\in \IC$ for each node, the Nakajima
variety is the symplectic reduction 
$$NQV(\Ga,V,\la)=\Rep(\Ga,V)\spqa{\la}\IG$$ 

\noindent
where $\la_i$ is identified with 
a scalar matrix in the centre of $\End(V_i)$, and so $\la$ is a central value of the moment map. 
Up to isomorphism the resulting symplectic variety is independent of the chosen orientation of $\Ga$.

Now suppose $G=\GL(V)$ for some vector space $V$ and we consider the reduction
$$\cB = \gaht\spq G$$
of some fission space (with $r=2$) with respect to $G$ (at the identity value of the $G$ moment map).
The space $\cB$ looks as in \eqref{eq: oppintn}, but with $T$ replaced by $H$, and $U_\pm$ by the unipotent radicals of a pair of opposite parabolics $P_\pm = H\cdot U_\pm\subset G$.
In particular $H=\prod_I\GL(V_i)$ for some grading $V=\bigoplus_I V_i$ of $V$.
Also, as a space, observe that $U_+\times U_-$ is isomorphic to $\Rep(\Ga,V)$ where $\Ga$ 
is the complete graph with nodes $I$. 
Further (as in \eqref{eq: oppintn}) $\cB$ is an open subset of $U_+\times U_-$.
In other words we have:

\begin{thm}
If $G=\GL(V), H=\prod_I\GL(V_i)$ then the space $\cB = \gaht\spq G$
is a quasi-Hamiltonian $H$-space and
is an open subset, the subset of ``invertible representations'' 
$$\Rep^*(\Ga,V)\subset \Rep(\Ga,V)$$
of the space of representations of the complete graph $\Ga$ on the vector space $V=\bigoplus_I V_i.$ 
\end{thm}

For example if $\Ga$ is a triangle and each $V_i$ is one dimensional, the reductions 
$\Rep^*(\Ga,V)\spq_q\,\, H$  (at values $q\in (\IC^*)^I$ of the moment map) are complex surfaces (real dimension four), 
and are multiplicative analogues of the corresponding Nakajima varieties, which in this case are the $A_2$ asymptotically locally Euclidean gravitational instantons of Gibbons--Hawking \cite{gib-haw78} (deformations of the minimal resolution of $\IC^2/\IZ_3$).
In fact these multiplicative spaces $\Rep^*(\Ga,V)\spq_q \,\, H$
are themselves gravitational instantons (complete \hk manifolds), due to \cite{wnabh}, but with different asymptotics at infinity (cf. the  conjectural list of the complex two-dimensional moduli spaces in \cite{ihptalk} \S3.2).

More generally one can consider also a pair of opposite parabolics $Q_\pm \subset H$
in the group $H$,
 and thus the fission space $\!\pap{H}{K}=\!\papk{H}{K}{1}$ determined by them 
 (where $K=Q_+\cap Q_-\subset H$ is the common Levi subgroup).
Then one can consider the gluing
$$\gaht \,\,\glue{H} \pap{H}{K}\quad = \quad (\!\gaht\,\, \fusion{H}\pap{H}{K})\spq H$$

\noindent
and in turn its reduction $(\!\!\gaht \,\glue{}\! \pap{H}{K})\spq G$, which is a quasi-Hamiltonian $K$-space.
Similarly to above, this reduction is naturally an open subset 
of the space of representations of a certain complete $k$-partite graph, and all complete $k$-partite graphs arise in this way ($k$ is the number of factors occurring in the intermediate group  $H$ here).
In turn, one may glue such graphs together to yield canonical open subsets of spaces of representations of general graphs, and then form the reduction to obtain many symplectic manifolds,  
and thus a full theory of ``multiplicative quiver varieties'' (see \cite{cmqv} for more details).

In this way a large class of fission varieties are determined by coloured graphs (we colour the graph so each coloured piece is a complete $k$-partite graph).
However only for special coloured graphs, such as the ``supernova graphs'', are these fission varieties actually wild character varieties, as shown in \cite{cmqv}. 
This leads to a theory of ``Dynkin diagrams'' for the wild character varieties (in fact for the 
full \hk manifolds with all their different structures).
One consequence of this theory, parameterising wild character varieties by certain graphs, is that 
the Weyl group and root system of the Kac--Moody algebra whose Dynkin diagram is the graph, plays an important  role  (see \cite{slims, cmqv}).

\subsection{Diagonal parts, and their multiplicative version}\label{ssn: diag parts}

Let $\cO\subset \g\cong \g^*$ be a (co)adjoint orbit for the complex reductive group $G$.
As is well-known $\cO$ is naturally a Hamiltonian $G$ space with moment map given by the inclusion 
$$ \cO \hookrightarrow \g^*.$$
The orbit method suggests we compare coadjoint orbits with representations $V$ of $G$, and given a representation 
one usually studies it by restricting to a maximal torus $T\subset G$ and decomposing it into simultaneous eigenspaces of $T$, i.e. the weight spaces of the representation.
In the classical/geometric version this means considering $\cO$ as a $T$-space, rather than as a $G$-space.
The action of $T$ is again Hamiltonian with moment map given by the dual of the inclusion $\lt\subset \g$:
$$\cO\subset \g^*\to \lt^*.$$
Using an invariant inner product to identify $\g^*\cong \g$ and $\lt^*\cong \lt$ the $T$ moment map becomes the map taking the ``diagonal part'' 
$$\de:\cO\subset \g\to \lt.$$
In the case of compact groups the image of such orbits under the diagonal part map was much studied 
by Kostant \cite{Kost73} and Heckman \cite{heckman-projs}, 
and this led to the Atiyah--Guillemin--Sternberg convexity theorem.

The rough idea is that when we view $V$ as a quantisation of $\cO$, the individual weight spaces 
$V_\la\subset V$ should be quantisations of the symplectic reductions 
\beq\label{eq: wt var}
\cO\spqa{\la} T\eeq
of $\cO$ by $T$ at the value $\la\in \lt^*$ of the moment map.
This is the viewpoint taken in Knutson's thesis \cite{knutson-these} who studied such ``weight varieties'' \eqref{eq: wt var} in the case of compact groups.

Now the multiplicative analogue of a coadjoint orbit is a conjugacy class $\cC\subset G$, 
and, by \cite{AMM}, it is a quasi-Hamiltonian $G$-space with moment map given by the inclusion.
However, we do not have a ``multiplicative diagonal part'' map $G\to T$ so it is not immediately clear how to convert a quasi-Hamiltonian $G$-space into a quasi-Hamiltonian $T$-space.

A solution to this is suggested by the fission operation. 
Let $\gat=\gato$ be a fission space (with $r=1$) associated to a pair of opposite Borel subgroups 
$B_\pm=T\cdot U_\pm \subset G$.
Then if we have a quasi-Hamiltonian $G$-space $M$ 
we may glue on the fission space $\gat$ to obtain a 
quasi-Hamiltonian $T$-space:
$$ M\quad \mapsto \quad 
M\glue{G} \gat = (M\fusion{G} \gat)\spq G.$$
Thus fission gives a way to break the structure group 
from $G$ to the maximal torus $T$ (hence the name ``fission'').
In the case $M=\cC\subset G$ the resulting quasi-Hamiltonian $T$-space 
$\cC\!\glue{}\!\!\! \gat$  may be identified with the intersection of $\cC$ with the 
``big cell'' $U_+TU_-\subset G$ and the moment map is then given by the (inverse of the) $T$-component.
(Of course since $T$ is abelian a quasi-Hamiltonian $T$-space is in particular a complex symplectic manifold.)

Of course it is not immediately clear that this is the right thing to do---it is strange to have to restrict to an open subset, since we didn't need to in the additive case.

However in the case $G=\GL_n(\IC)$ the weight varieties have an alternative description, and it is clear what the multiplicative versions of them are, and they are isomorphic to the ``multiplicative weight varieties'' 
$(\cC\,\glue{}\gat) \spq_q \,\,T$:

\begin{thm}\label{thm: add duality}
For any weight variety $\cO\spq_{\la}\,\, T$ for the group 
$G=\GL_n(\IC)$ there is an integer $m$, and 
coadjoint orbits $\cO_1,\ldots,\cO_n, \cO_\infty\subset \gl_m(\IC)^*$
of $\GL_m(\IC)$ such that $\cO\spq_\la\,\, T$ is isomorphic to the 
variety
\beq\label{eq: dwv}
(\cO_1\times\cdots \times\cO_n)\spqa{\cO_\infty} \GL_m(\IC).
\eeq
\end{thm}

More precisely this holds at the level of subsets of {\em stable} points, and in the case at hand each of the orbits $\cO_i$ can be taken to be an orbit of rank one 
matrices ($1\le i\le n$).
This arises by thinking about the reduction of $T^*\Hom(\IC^m,\IC^n)$ via the natural action of $T\times \GL_m(\IC)$, in two different ways.
Said differently these are two of the ways to ``read'' a star-shaped Nakajima quiver variety in terms of moduli spaces of connections (cf. \cite{iastalk07-nom, rsode} and  the more general results in \cite{slims} \S9)---in fact this
type of result comes about by thinking about the Fourier--Laplace transform 
on spaces of meromorphic connections on the Riemann sphere---see \cite{nlsl} Diagram 1 and references therein---in other terms this isomorphism is also a precise version of Harnad's duality \cite{Harn94}.

In particular, from the moduli space interpretation of the weight variety and its dual \eqref{eq: dwv}, we know immediately what their ``multiplicative versions'' are, namely we look at the corresponding (wild) character varieties, on the other side of the (irregular) Riemann--Hilbert correspondence:

First \eqref{eq: dwv} is a moduli space of logarithmic connections on the trivial rank $m$ vector bundle on the Riemann sphere, having a first order pole at $n$ distinct points $a_1,\ldots,a_n\in\IC$, and a further pole at $\infty$ 
(cf. e.g. \cite{smid} \S2).
The corresponding irregular curve is just the Riemann sphere with marked points $a_1,\ldots,a_n,\infty$ (with trivial irregular types), and so the corresponding (tame) character varieties are of the form 
\beq\label{eq: tcv}
(\cC_1\fus{}\cdots \fus{} \cC_n)\spqa{\cC_\infty} \GL_m(\IC)
\eeq

\noindent
for conjugacy classes $\cC_1,\ldots,\cC_n,\cC_\infty \subset \GL_m(\IC)$.
These are the multiplicative versions of \eqref{eq: dwv}.
For sufficiently generic (nonresonant) orbits  the Riemann--Hilbert map is a
holomorphic map from \eqref{eq: dwv} to \eqref{eq: tcv}.

Secondly the complex weight varieties $\cO\spq_{\la}\,\, T$ are moduli spaces of 
meromorphic connections on the trivial bundle over the Riemann sphere, with poles only at $0,\infty$, having trivial irregular type at $\infty$, and irregular type 
$Q=A/z$ at $z=0$, for some regular semi-simple element $A\in \lt$.
The orbit $\cO$ specifies the orbit of the residue at $\infty$, and the weight 
$\la\in \lt\cong \lt^*$ specifies the residue of the formal diagonalisation 
 $dQ+\la dz/z$
of
\eqref{eq: fnf}
 at $z=0$ (viewing $\{\la\}$ as a coadjoint orbit of $T$).
It is a nice exercise to see that the moduli space of such connections is the symplectic quotient $\cO\spq_{\la} \,\,T$ (i.e. that fixing the residue of the normal form at $0$ amounts to fixing the diagonal part, cf. \cite{smapg}---note that here we are working with singular connections on bundles on the compact curve $\IP^1$, rather than algebraic connections on the punctured curve).  
The corresponding irregular curve is thus $(\IP^1, (0,\infty),(Q,0))$ with two marked points and the given irregular types.
Moreover, unwinding the definitions,
the corresponding wild character varieties are of the form
$$(\cC \glue{G} \gat)\spqa{q} T$$

\noindent
where $q\in T$ is a point (i.e. a conjugacy class of $T$). In other words they are 
the varietes we proposed above to call the ``multiplicative weight varieties''.
Again for sufficiently generic (nonresonant) orbits the irregular Riemann--Hilbert map is a holomorphic map from the weight variety to its multiplicative version (in fact this holds for all reductive groups---and so gives a good justification that we have the right definition).
A key point is that the irregular type $Q$ thus provides a natural mechanism to break the structure group from $G$ to $T$.

A further justification arises since the multiplicative version of Theorem \ref{thm: add duality} holds:

\begin{thm}\label{thm: mult duality}
For any multiplicative weight variety $(\cC\glue{G}\!\!\gat)\spq_{q}\,\, T$ for the group 
$G=\GL_n(\IC)$ there is an integer $m$, and 
conjugacy classes $\cC_1,\ldots,\cC_n, \cC_\infty\subset \GL_m(\IC)$
of $\GL_m(\IC)$ such that $(\cC\glue{G} \gat)\spq_{q} \,\,T$ is isomorphic to the 
variety
\beq\label{eq: mdwv}
(\cC_1\fus{}\cdots \fus{}\cC_n)\spqa{\cC_\infty} \GL_m(\IC).
\eeq
\end{thm}

Strictly speaking, again we should add the caveat of restricting to stable points here (for sufficiently generic conjugacy classes everything is stable, cf. \cite{gbs} \S9.2). 
Such isomorphisms arise  by computing the action of the Fourier--Laplace transform on Betti data, and many generalisations are possible (see \cite{cmqv}). 
In other words the fission space $\gat$ plays a similar role 
to that of the space 
$T^*\Hom(\IC^m,\IC^n)$ in the additive version above.

\begin{rmk}
As already noted in \cite{iastalk07-nom, rsode} the additive duality above 
may be viewed as a generalisation of (the complexification of) 
the duality used by Gelfand--MacPherson \cite{GelfMacP} in their study of Grassmannians
and the dilogarithm.
In fact in \cite{hausmann-knutson-poly} (7.5)  Hausmann--Knutson enquired about the generalisation of the 
Gelfand--MacPherson duality to the multiplicative case, for compact groups---whilst this is not clear, we are saying that the complexification has a nice multiplicative version, and a longer history
(going back to Fourier--Laplace, or really to Poincar\'e and Birkhoff for the study of the Fourier--Laplace transform in the context of irregular connections---see the references in  \cite{BJL81}).
Such ideas were used to give the first intrinsic derivation of the Regge symmetry of the 6j-symbols in 
\cite{roks}.
\end{rmk}

\begin{rmk}
Again in the context of compact groups, the multiplicative version of the diagonal part $\de$ map, is the ``Iwasawa projection'' $\wh\de$, which appears in Kostant's nonlinear convexity theorem \cite{Kost73}.
The possibility of finding a map intertwining $\de$ and $\wh \de$ (twisting the exponential map appropriately) was discovered by Duistermaat \cite{Dui84}.
From the moduli space viewpoint, such twistings arise naturally as certain connection matrices between the two poles of the irregular connections \cite{smapg} (see also the review in \cite{hdr}).
In a similar spirit the irregular Riemann--Hilbert map gives many {\em naturally occuring} Ginzburg--Weinstein isomorphisms $\k^*\cong K^*$ (see op. cit.).
\end{rmk}

Many generalisations of this set-up can be (and have been) considered. 
Some of the simplest possibilities are as follows.

1) The condition that the leading coefficient $A\in\lt$ 
(appearing in the irregular type $Q$) 
is regular implies that the centralizer group $H$ is a maximal torus 
and on the other side of the duality implies that the first $n$ orbits/conjugacy classes are of minimal dimension.
If we relax this condition, allowing any $A\in\lt$, the duality still holds, and this leads to the  fact that 
{\em any} genus zero tame character variety (for a general linear group)
is {\em isomorphic} 
to a wild character variety of the form
$$(\cC\glue{G}\gah)\spqa{\breve \cC} H\quad\cong\quad
 \cC\glue{G}\gah \glue{H} \breve\cC'$$
for some classes $\cC\subset G, \breve\cC\subset H$ (where $\breve\cC'$ denotes the inverse conjugacy class to $\breve \cC$),
i.e. it is isomorphic to a symplectic wild character variety for an irregular curve 
of the form $\Si=(\IP^1,(0,\infty),(Q,0))$ with $Q=A/z$ for some $A\in \lt$
(see \cite{cmqv}  Theorem 1.5 for more details---the additive version is a special case of results of \cite{slims} \S9).

Thus in terms of pole orders, in genus zero for general linear groups, we have an equivalence
\beq\label{eq: full duality1}
1+1+\cdots+1\qquad \Longleftrightarrow \qquad 2+1
\eeq
between the tame case and the simplest irregular case.
Even though the Stokes data might look strange (involving unipotent subgroups), we actually get moduli spaces isomorphic to well-known/much studied spaces. 
Of course there are many wild character varieties beyond this simple case, that are {\em not} isomorphic to a tame case.

2) The next simplest case is, symbolically (in terms of pole orders):
$$2+1+\cdots+1\qquad \Longleftrightarrow \qquad 2+1+\cdots +1$$ 
i.e. allowing a pole of order two 
(from a nontrivial irregular type of the form $A/z$ at one point) on both sides, and allowing arbitrarily many simple poles on each side.
The additive version of this duality is the context of 
Harnad's duality \cite{Harn94}.
Here the additive moduli spaces are of the form
$$(\cO_1\times \cdots\times \cO_n)\spqa{\breve \cO} H$$

\noindent
for coadjoint orbits $\cO_i$ of $G=\GL_m(\IC)$, and $\breve \cO$ of $H\subset G$.
The duality changes the rank $m$ and 
exchanges the number of orbits $n$ with the number of factors 
of the groups $H$ (i.e. the number of distinct eigenvalues of $A$), and has a nice interpretation in terms of quivers (cf. \cite{iastalk07-nom, rsode, slims}).
The corresponding multiplicative moduli spaces (wild character varieties) 
are of the form
$$(\cC_1\fusion{G} \cdots\fusion{G} \cC_n\glue{G}\gah)\spqa{\breve \cC} H$$

\noindent
for conjugacy classes $\cC_i$ of $G=\GL_m(\IC)$, and $\breve \cC$ of $H\subset G$.
The multiplicative duality (see \cite{cmqv}) again has a nice interpretation in terms of quivers, and comes from the computation of the action of the Fourier--Laplace transform on Betti data.

This example shows in fact that the eigenvalues of $A$ (the coefficient of the irregular type $Q=A/z$) correspond on the other side to the 
positions $a_1,\ldots,a_n\in \IC$ of the poles (of the tame singularities in the complex plane). This provided motivation for the notion of ``irregular curve'', since the moduli of the pointed curve $(\IP^1,(a_1,\ldots,a_n,\infty))$ 
matches with the choice of the irregular type $Q$ on the other side (so one gets a precise correspondence between the moduli of {\em irregular curves} on both sides, including the pole positions and the irregular type on both sides).
Said differently, in these examples 
there are two groups of pole positions, only one of which is expressed as actual pole positions at any given time (one passes to the other side of the duality to realize the irregular type as genuine pole positions). 
Thus in this example it is not at all surprising that the irregular type 
``behaves the same way'' as the moduli of the pointed curve, since under the duality 
they are swapped over.
Thus the notion of irregular curve can be seen as a natural generalisation of the notion of curve with marked points.
This viewpoint  is especially useful when we consider deformations/braiding below.

3) The next simplest class of examples involves connections having a pole of order three, i.e. an irregular type of the form $Q=A/z^2+B/z$ for some $A,B\in \lt\subset \gl_n(\IC)$.
The simplest set-up is to take $\Si=(\IP^1,0,Q)$, i.e. the Riemann sphere with just one marked point, with such $Q$.
We will call such a curve $\Si$ a ``type $3$ irregular curve''.
The class of moduli spaces which appear in this case is incredibly rich, since for example:

\begin{thm}[\cite{cmqv}]
All of the (stable) wild character varieties which occur in cases 1) or 2) above, i.e. with poles of type 
$2+1+\cdots+1$,
are isomorphic to a wild character variety of a type $3$ irregular curve $\Si$.
\end{thm}

In brief whereas originally we had two groups of ``pole positions'' which could be swapped over, 
we now have $k$ groups of ``pole positions'', where $k$ is the number of distinct eigenvalues of
the leading term $A$ of the irregular type $Q=A/z^2+B/z$.
Any one of these groups of ``pole positions'' can be realised as an actual tuple of pole positions of  connections with poles of the form $3+1+\cdots+1$ (on a lower rank bundle, and with one less eigenvalue in the leading term at the pole of order $3$).
The full story (see \cite{slims}) is that in fact we may act on the eigenvalues of $A$ with 
arbitrary M\"obius transformations, moving them around in the Riemann sphere, and when a particular eigenvalue is moved to infinity then the corresponding group of ``pole positions'' 
are realised as actual pole positions.
In the case $k=2$ we get moduli spaces isomorphic to those considered earlier, but in general we get new moduli spaces (with more collections of pole positions).
For more details see  \cite{cmqv},  or  \cite{slims} which studied the simpler additive version 
of this set-up in detail, and related such isomorphisms between moduli spaces to automorphisms of the Weyl algebra.

\section{Wild mapping class groups}

The usual story of braid/mapping class group actions on (tame) character varieties may be described as follows.

Suppose $\pi:\Si\to \IB$ is a family of smooth compact Riemann surfaces over a base $\IB$, 
so that each fibre $\Si_b=\pi^{-1}(b)$ is a smooth compact Riemann surface.
Suppose further each fibre has $m$ distinct marked points, and these points vary smoothly, so that we have sections $a_i:\IB\to \Si$ for $i=1,\ldots,m$, and the subvarieties $a_i(\IB)\subset \Si$ do not meet 
each other.

Then, given a complex reductive group  $G$, and a point $b\in \IB$ we get a character variety 
$$\MB(b) = \Hom(\pi_1(\Si^\circ_b),G)/G$$
which has a natural Poisson structure (given the choice of invariant bilinear form on $\g$),
where $\Si^\circ_b=\Si_b\setminus\{a_1(b),\ldots,a_m(b)\}$ is the $m$-punctured surface over $b\in \IB$.

These character varieties $\MB(b)$ then fit together into the fibres of a fibre bundle
$$\pr:\MB \to \IB$$
with $\pr^{-1}(b)=\MB(b)$
with a natural flat (nonlinear/Ehresmann) connection on it, preserving the Poisson structure.
One way to see this is to note that locally, over a small open ball $U$ in $\IB$, we may identify 
the fibres $\MB(b_1)\cong \MB(b_2)$ for $b_1,b_2\in \U$, by using the {\em same} loops generating 
$\pi_1(\Si^\circ_{b_i})$, and thereby identifying monodromy representations.
Said differently if we write $\Si^\circ= \Si \setminus \bigcup_1^m a_i(\IB)$
for the family of punctured curves,
then 
the two inclusions 
$\Si^\circ_{b_i} \hookrightarrow \Si^\circ\bigl\vert_U$
are homotopy equivalences 
and so the fundamental groups (and thus the character varieties) are identified:
$$\pi_1(\Si^\circ_{b_1})\quad  \cong\quad \pi_1(\Si^\circ\bigl\vert_U)\quad \cong\quad \pi_1(\Si^\circ_{b_2}).$$
Globally this gives a presentation of the bundle $\pr:\MB \to \IB$ in terms of local trivialisations with constant, Poisson, clutching maps, and thus  the bundle comes with a flat nonlinear Poisson connection.

Integrating this connection gives a Poisson action of the fundamental group $\pi_1(\IB,b)$ 
of the base $\IB$,  on the fibre $\MB(b)$ over the base point $b\in \IB$.  

For example taking $\IB=\IC^m\setminus \text{diagonals}$ to be the configuration space of $m$-tuples of points in the plane leads to a Poisson action of the pure Artin braid group $\pi_1(\IB)$ on the genus zero character varieties, and in the higher genus case one obtains action of the mapping class group.
More generally one can consider unordered points (replacing the sections $a_i$ by a multisection) 
and obtain actions of the full (non-pure) Artin/type $A$  braid groups.

For us the key point is that a similar story also holds true in the more general context of irregular curves, not just the tame case considered above, and there is interesting braiding of the irregular types.
The main statement is as follows.

\begin{thm}(\cite{gbs})
The fundamental group of any space $\IB$ of admissible deformations of an irregular curve $\Si_b$
acts on the wild character variety $\MB(\Si_b)$ of $\Si_b$ by algebraic Poisson automorphisms.
\end{thm}

More precisely (see \cite{gbs}) the wild 
character varieties again fit together into the fibres of 
fibre bundle over $\IB$, with a natural complete flat nonlinear Poisson connection, 
and the action of $\pi_1(\IB)$ arises from integrating this connection.
An earlier analytic approach to this result appeared in \cite{thesis, smid}, and that was given a supersymmetric interpretation 
in \cite{witten-wild}.

It is clear that such a result will not hold under arbitrary deformations of the irregular types $Q_i$, 
since the dimensions of the wild character varieties will change 
(cf. \eqref{eq: big fusion prod} and \cite{gbs} Rmk 9.12 for the formula for the dimension).
The condition of only allowing ``admissible'' deformations of the irregular curve generalises the notion of deforming a curve with marked points such that the curve remains smooth and the marked points do no coalesce. It brings into play certain hyperplanes in the space of deformations of irregular types, and thus leads to deformation spaces with interesting fundamental groups.
The general definition is as follows.

\begin{defn}
An ``admissible deformation'' of an irregular curve $\Si_0$, consists of a family of irregular curves (with one fibre isomorphic to $\Si_0$),  such that 1) each fibre $\Si_b$  is smooth, 2) the marked points remain distinct, and 3) for any $i=1,\ldots,m$ and  any 
root $\al\in \cR\subset \lt^*$ of $\g$, the order of the pole of
$$\al \circ Q_i$$
does not change, where $Q_i$ is the irregular type at the $i$th marked point.
\end{defn}

Simple examples of admissible deformations were considered by Jimbo--Miwa--Ueno 
\cite{JMU81} who looked at the genus zero case for general linear groups, and imposed the condition that the leading coefficient at each irregular singularity had distinct eigenvalues.
Since for general linear groups the roots correspond to the differences of the eigenvalues, this is a special type of admissible deformation when  all the pole orders $\al \circ Q_i$ are the same (at each marked point).

For example taking the case of poles of type $2+1$ (as in \S\ref{ssn: diag parts})
for the group $G=\GL_n(\IC)$ we could take 
$$\IB=\lt_\reg=\IC^n\setminus\text{diagonals}$$
so an element $A\in \IB$ determines an irregular curve 
$\Si_A = (\IP^1,(0,\infty),(A/z,0))$
and so we obtain an admissible family of irregular curves over $\IB$ in this way, and thus a Poisson action of the Artin braid group 
$$\pi_1(\IB) = \pi_1(\lt_\reg)$$
on the corresponding wild character varieties.
In this case (for general linear groups) the duality described in \S\ref{ssn: diag parts} is braid group equivariant: all these wild character varieties are isomorphic to tame genus zero character varieties and the braid group actions match up.

If we start instead with arbitrary $A\in \lt$ in this example,  
then, as in \eqref{eq: full duality1}, all the tame genus zero (general linear) character varieties are 
obtained as wild character varieties, and again the duality is braid group equivariant (the admissibility condition is that no further eigenvalues of $A$ should coalesce). 

In these examples only  Artin braid groups (type $A$) are obtained: on the tame side this is because we are moving around pole positions in the complex plane, and on the irregular side it is because we are working with the structure group $\GL_n(\IC)$.
Thus although it is not at all clear how to generalise the tame picture to obtain generalised braid groups (associated to other groups $G$), it is clear on the irregular side: we just need to work with irregular connections on principal $G$-bundles, rather than vector bundles, and extend the theory of Stokes data to this context.

This was carried out in \cite{bafi}---in the simplest case one considers 
$$\IB=\lt_\reg=\{A\in \lt \st \al(A)\neq 0 \text{ for all roots $\al$}\} $$
where $\lt\subset \g=\Lie(G)$ is a Cartan subalgebra of our fixed complex reductive group $G$ (for example $G=E_8\times G_2\times\SO_8$).
Thus an element $A\in \IB$ determines an irregular curve 
$\Si_A = (\IP^1,(0,\infty),(A/z,0))$
and so we obtain an admissible family of irregular curves in this way, and thus a Poisson action of the pure $\g$-braid group 
$$\pi_1(\IB) = \pi_1(\lt_\reg)=P_\g$$
on the corresponding wild character varieties.
(One can extend this to the full braid groups $B_\g=\pi_1(\lt_\reg/W)$ 
by considering ``bare irregular curves'', as discussed in \cite{gbs}.)

In this simplest case such a braid group action had been seen before. More precisely if we consider the local version with such irregular type $Q=A/z$ on a disk (around $z=0$) then the Poisson action lifts to the framed wild character varieties (adding in a $T$-framing at $0$), which are the Poisson manifolds 
$\wt G^*\to G^*\to G^\circ$.
Such a Poisson action (of $B_\g$ on $G^*$) was written down in 
\cite{DKP} by computing the quasi-classical limit of the so-called quantum Weyl group action (which is essentially a $B_\g$ action on the Drinfeld--Jimbo quantum group, constructed by generators and relations by Lusztig, Soibelman, Kirillov--Reshetikhin).

\begin{thm}(\cite{bafi})
The action of $B_\g$ on $G^*$ derived by DeConcini--Kac--Procesi \cite{DKP} coincides with the geometric action obtained by computing the monodromy action of the natural nonlinear connection on the bundle of (framed) wild character varieties. 
\end{thm}

Of course there are many examples beyond this case. 
For example in the notation of \S\ref{ssn: diag parts} the above examples are of type $2+1$, or  $1+\cdots +1=1^m$.
The next simplest $\GL_n$ case is of type $2+1^m$ and involves two groups of ``pole positions'', and thus one can consider admissible families over 
$$\IB = (\IC^n\setminus\diagonals)\times(\IC^m\setminus\diagonals).$$
Next, as mentioned in \S\ref{ssn: diag parts}, these are isomorphic 
to special cases of the wild character varieties of type $3$ irregular curves.
In general type $3$ irregular curves may have $k$ groups of ``pole positions'', and thus one may consider admissible families of type $3$ irregular curves over 
$$\IB = (\IC^{n_1}\setminus\diagonals)\times\cdots\times (\IC^{n_k}\setminus\diagonals).$$
for any integers $n_i$ (as in \cite{rsode, slims}).

In general, given any irregular curve $\Si$, one would 
consider a universal space $\gM(\Si)$
of admissible deformations of $\Si$, generalising the moduli space/stack of 
smooth curves with marked points in the tame case, 
and then obtain a flat Poisson bundle of wild character varieties over $\gM(\Si)$,
and thus a Poisson action of the fundamental group $\pi_1(\gM(\Si))$ (the ``wild mapping class group'') on the wild character variety $\MB(\Si)$.

\appendix
\section{Singular directions and Stokes groups}\label{apx: sd and sg}

Suppose we have a smooth curve $\Si$, a point $a\in \Si$ and an irregular type $Q$ at $a$.
This appendix will recall %
the definition of the singular directions
$\IA\subset \partial$ and the Stokes groups $\ISto_d(Q)\subset G$ for each $d\in \IA$.

Recall $T\subset G$ is a maximal torus, and let $\lt\subset \g$ denote the corresponding Lie algebras.
Let $\cR\subset \lt^*$ be the roots of $\g$ 
so that 
$$\g = \lt \oplus\bigoplus_{\al\in \cR} \g_\al$$
where $\g_\al = \{ X\in \g \st [Y,X] = \al(Y)X\text{ for all $Y\in \lt$}\}$ is the (one dimensional) root space of $\al\in \cR$.
For each root $\al \in \cR$  let
$$q_\al = \al \circ Q$$
denote the $\al$ component of $Q$.
Define the degree $\deg(q_\al)\in \IZ_{\ge 0}$ of 
$q_\al$ to be its pole order at $a$.
Let $\wt \Si\to\Si$ denote the real oriented blow-up at $a$, replacing $a$ with the circle $\partial$ of real oriented tangent directions at $a$.

\begin{defn} \label{def: sing dirn}
A direction $d\in \partial$ is a {\em ``singular direction supported by $\al$'',} (or an {\em ``anti-Stokes direction''})  if
$q_\al$ has a pole, and 
$\exp(q_\al(z))$ 
has maximal decay as $z\to 0$ in the direction $d$.
\end{defn}

Thus if $c_\al/z^k$ is the most singular term of $q_\al$, with $c_\al\in\IC^*$ (in a local coordinate $z$ vanishing at $a$), 
and $k\ge 1$, then the singular directions supported by $\al$
are those along which the function $c_\al/z^k$ is real and negative.
Let $\IA\subset \partial$ be the finite set of singular directions (for all roots $\al$).
If $d\in\IA$ is a singular direction
let $$\R(d)\subset \cR$$ denote the (nonempty) subset of roots supporting $d$.

\begin{defn}
The Stokes group $\ISto_d(Q)\subset G$ is the subgroup of $G$ generated by the root groups $\exp(\g_\al)$ for all $\al\in \cR(d)$.
\end{defn}

It follows that $\ISto_d(Q)$ is a unipotent subgroup of 
$G$, 
normalised by $H=C_G(Q)$, of dimension equal to $\#\cR(d)$.

For example (similarly to \cite{bafi}) if  
$\Si=\Delta$ is the unit disc in $\IC$, and $Q=-A/z$ (so that $dQ=Adz/z^2$) then $\IA$ is the set of rays from $0$ to the nonzero points  in the set
$\langle \cR , A \rangle \subset \IC$,
where the angled brackets %
denote the natural pairing between $\lt^*$ and $\lt$.
Further $\cR(d)$ is then the set of roots landing 
on the ray $d$ (and in $\IC^*$)
under the map $\langle \,\cdot\,, A \rangle:\cR\to \IC$.
If we replace $Q=-A/z$ by $Q'=-A/z^k$ in this example, the picture is similar except each root is repeated $k$ times as we go around the circle (i.e. we just pull-back under the map $z\mapsto z^k$).
The cases where $Q$ has multiple terms 
(so that the components $q_\al$ may have different degrees as $\al$ varies) are more complicated since different roots repeat with different frequencies as we go around the circle. Nonetheless the resulting spaces of Stokes data can be obtained by gluing together the simpler spaces with just one term (cf. \cite{gbs} Prop. 7.12).

The set of ``levels'' of an irregular type $Q$ is the set of degrees  
$\{\deg(q_\al)\st \al \in \cR\}\setminus \{0\}\subset \IZ_{>0}$
that occur.
Thus the level of a root $\al$ is the number of times it is repeated as we go around the circle $\partial$.
By gluing the fission spaces end to end, an induction is possible with respect to the number of levels.
Note that if the leading term $A_r$ of $Q$ is regular semisimple 
then $Q$ has just one level: this condition means that 
$\al(A_r)\neq 0$ for all roots. Thus this condition is a simple way to avoid the complications that arise in the multi-level cases.

\section{}\label{apx: morecan}
Recall the notion of Stokes local system involved a local system on a surface 
$\wt \Si\subset \wh \Si\to \Si$, where $\wh\Si$ was the real oriented blow-up of $\Si$ at each of the marked points $a_i$.
A slightly different approach to defining the surface  $\wt \Si$ is as follows. (It is clearly homotopy equivalent to the approach above, but has the benefit of not having to worry about specifying exactly where to place the punctures $e(d)$.)
Given $\Si$ and the real blow-up $\wh \Si\to \Si$ with boundary circles $\partial_1,\ldots,\partial_m$ as above, glue on an annulus $\IH_i=[1/2,1]\times S^1$ to each boundary circle $\partial_i$. 
This yields a new surface $\wh\Si'\to \Si$ over $\Si$ (in effect replacing $a_i\in \Si$ by $\IH_i$) so that $\partial_i$  is now in the interior of 
$\wh\Si'$. Now define $\wt \Si$ by puncturing 
$\wh\Si'$ at each singular direction $d\in \IA_i\subset \partial_i$ for each $i$.
The picture is as in Figure \ref{fig: halo} except now the dashed circle is $\partial_1$, %
the black dots are deleted and the singular directions are the white dots, where we puncture.

\renewcommand{\baselinestretch}{1}              %
\normalsize
\bibliographystyle{amsalpha}    \label{biby}
\bibliography{../thesis/syr} 
{\small
\noindent
Institut des Hautes \'Etudes Scientifiques et CNRS,

\noindent
Le Bois-Marie, 35 route de Chartres, 
91440 Bures-sur-Yvette, France

\noindent
boalch@ihes.fr

\end{document}